\documentclass[final]{siamltex} 	
\usepackage{amsmath,amssymb,mathdots,extarrows,cite,mathtools}
\usepackage{graphicx}
\usepackage{tikz-cd}
\usepackage{algorithm}
\usepackage{algpseudocode}

\newcommand{\mzero}{\mathbf O }
\newcommand{\vzero}{\mathbf 0 }

\newtheorem{remark}[theorem]{Remark}
\newtheorem{example}[theorem]{Example}

\newcommand{\hatdSmR}[2]{\widehat{\mathbf d}^{#1}}

\newcommand{\Diag}[1]{\textup{\text{Diag}}(\mathbf #1)}
\newcommand{\Matr}[1]{\textup{\text{Matr}}( #1)}
\newcommand{\UnMatr}[1]{\textup{\text{Tens}}( #1)}
\newcommand{\vect}[1]{\textup{\text{vec}}( #1)}

\newcommand{\rdm}[1]{{\mathcal F}_{#1}}
\newcommand\Algphase[1]{%
\vspace*{-.4\baselineskip}\Statex\hspace*{\dimexpr-\algorithmicindent-2pt\relax}\rule{\textwidth}{0.4pt}%
\Statex\hspace*{-\algorithmicindent}{#1}%
\vspace*{-.7\baselineskip}\Statex\hspace*{\dimexpr-\algorithmicindent-2pt\relax}\rule{\textwidth}{0.4pt}%
}
\newcommand{\setalglineno}[1]{%
  \setcounter{ALG@line}{\numexpr#1-1}}

\newcommand{\RTm}[1]{{\mathbf R}_{#1}(\mathcal T)}
\newcommand{\TTm}[1]{{\mathbf Q}_{#1}(\mathcal T)}

\renewcommand{\algorithmicrequire}{\textbf{Input:}}
\renewcommand{\algorithmicensure}{\textbf{Output:}}
\renewcommand{\algorithmicrepeat}{}
\renewcommand{\algorithmicuntil}{}
\makeatletter
\newif\if@borderstar
\def\bordermatrix{\@ifnextchar*{%
\@borderstartrue\@bordermatrix@i}{\@borderstarfalse\@bordermatrix@i*}%
}
 \def\@bordermatrix@i*{\@ifnextchar[{\@bordermatrix@ii}{\@bordermatrix@ii[()
]}}
 \def\@bordermatrix@ii[#1]#2{%
 \begingroup
 \m@th\@tempdima8.75\p@\setbox\z@\vbox{%
 \def\cr{\crcr\noalign{\kern 2\p@\global\let\cr\endline }}%
 \ialign {$##$\hfil\kern 2\p@\kern\@tempdima & \thinspace %
 \hfil $##$\hfil && \quad\hfil $##$\hfil\crcr\omit\strut %
 \hfil\crcr\noalign{\kern -\baselineskip}#2\crcr\omit %
 \strut\cr}}%
 \setbox\tw@\vbox{\unvcopy\z@\global\setbox\@ne\lastbox}%
 \setbox\tw@\hbox{\unhbox\@ne\unskip\global\setbox\@ne\lastbox}%
 \setbox\tw@\hbox{%
 $\kern\wd\@ne\kern -\@tempdima\left\@firstoftwo#1%
 \if@borderstar\kern2pt\else\kern -\wd\@ne\fi%
 \global\setbox\@ne\vbox{\box\@ne\if@borderstar\else\kern 2\p@\fi}%
 \vcenter{\if@borderstar\else\kern -\ht\@ne\fi%
 \unvbox\z@\kern-\if@borderstar2\fi\baselineskip}%
 \if@borderstar\kern-2\@tempdima\kern2\p@\else\,\fi\right\@secondoftwo
#1 $%
 }\null \;\vbox{\kern\ht\@ne\box\tw@}%
 \endgroup
 }
 \makeatother

\title{Canonical polyadic decomposition of third-order tensors: reduction to generalized eigenvalue decomposition
\thanks{Research supported by: (1) Research Council KU Leuven:
GOA-MaNet,  CoE EF/05/006 Optimization in Engineering (OPTEC),
CIF1,  STRT 1/08/23,  (2) F.W.O.: projects G.0427.10N, G.0830.14N, G.0881.14N, (3) the Belgian Federal Science Policy
Office: IUAP P7/19 (DYSCO,  ``Dynamical systems,  control and
optimization'',  2012--2017), (4) DBOF/10/015.}}

\author{Ignat Domanov\footnotemark[2] \footnotemark[3] \footnotemark[4]
\and
Lieven De Lathauwer\footnotemark[2] \footnotemark[3] \footnotemark[4]}
\begin{document}
\maketitle
\renewcommand{\thefootnote}{\fnsymbol{footnote}}

\footnotetext[2]{Group Science, Engineering and Technology, KU Leuven - Kulak,
E. Sabbelaan 53, 8500 Kortrijk, Belgium
 ({\tt ignat.domanov,\ lieven.delathauwer@kuleuven-kulak.be}).}
\footnotetext[3]{Dept. of Electrical Engineering  ESAT/STADIUS KU Leuven,
Kasteelpark Arenberg 10, bus 2446, B-3001 Leuven-Heverlee, Belgium
({\tt lieven.delathauwer@esat.kuleuven.be}).}
\footnotetext[4] {iMinds Future Health Department.}

\begin{abstract}
Canonical Polyadic Decomposition (CPD) of a third-order tensor is  decomposition in a minimal number of rank-$1$ tensors.
We call an algorithm algebraic if it is guaranteed to find the decomposition when it is exact and if it only relies
on standard linear algebra (essentially sets of linear equations and matrix factorizations).
The known algebraic algorithms for the  computation of the CPD are limited to cases where at least one of the  factor matrices  has full column rank.
In the paper we present an algebraic algorithm for the computation of the CPD in cases where none of the factor matrices has full column rank.
In particular, we show that if the famous Kruskal condition holds, then the CPD can be found algebraically.
\end{abstract}

\begin{keywords}
Canonical Polyadic Decomposition, Candecomp/Parafac Decomposition, tensor, Khatri-Rao product, compound matrix, permanent, mixed discriminant
\end{keywords}

\begin{AMS}
15A69, 15A23
\end{AMS}

\pagestyle{myheadings}
\thispagestyle{plain}
\markboth{IGNAT DOMANOV AND LIEVEN DE LATHAUWER}{Canonical polyadic decomposition: reduction to GEVD}

\section{Introduction}
\subsection{Basic notations and terminology}\label{Subsection1.1}
Throughout the paper $\mathbb R$ denotes the field of real numbers and $\mathcal T=(t_{ijk})\in\mathbb R^{I\times J\times K}$ denotes a third-order tensor with frontal slices $\mathbf T_1,\dots,\mathbf T_K\in\mathbb R^{I\times J}$;
$r_{\mathbf A}$, $\textup{range}(\mathbf A)$, and
 $\textup{ker}(\mathbf A)$ denote the rank, the range, and the null space of a matrix $\mathbf A$, respectively;
$k_{\mathbf A}$ (the $k$-rank  of  $\mathbf A$) is the largest number  such that
every subset of $k_{\mathbf A}$ columns of the matrix $\mathbf A$ is  linearly independent;
$\omega(\mathbf d)$ denotes the number of nonzero entries of a vector $\mathbf d$;
$\text{span}\{\mathbf f_1, \dots, \mathbf f_k\}$ denotes the linear span of the vectors $\mathbf f_1, \dots, \mathbf f_k$;
$\mzero_{m\times n}$, $\vzero_m$, and $\mathbf I_n$ are the zero $m\times n$ matrix, the zero $m\times 1$ vector, and the $n\times n$ identity matrix, respectively;
$C_n^k$ denotes the binomial coefficient,  $C_n^k=\frac{n!}{k!(n-k)!}$;
$\mathcal C_m(\mathbf A)$ (the $m$-th compound matrix  of $\mathbf A$) is  the  matrix containing the determinants of all $m\times m$ submatrices of $\mathbf A$, arranged with the submatrix index sets in lexicographic order (see \S \ref{sec:compandpermcomp} for details).

The outer product $\mathbf a\circ\mathbf b\circ\mathbf c\in\mathbb R^{I\times J\times K}$ of three nonzero vectors $\mathbf a\in\mathbb R^I$,
$\mathbf b\in\mathbb R^J$ and $\mathbf c\in\mathbb R^K$ is
 called {\em rank-$1$} tensor ($(\mathbf a\circ\mathbf b\circ\mathbf c)_{ijk}: = a_i b_j c_k$ for all values of the indices).

 A {\em Polyadic  Decomposition} of  $\mathcal T$ expresses $\mathcal T$ as a  sum of rank-$1$ terms:
\begin{equation}
\mathcal T=\sum\limits_{r=1}^R\mathbf a_r\circ \mathbf b_r\circ \mathbf c_r,
 \label{eqintro2}
\end{equation}
where $\mathbf a_r \in \mathbb R^{I}$, $\mathbf b_r \in \mathbb R^{J}$, $\mathbf c_r \in \mathbb R^{K}$, $1 \leq r \leq R$.
If the number $R$ of rank-1 terms in \eqref{eqintro2} is minimal, then
\eqref{eqintro2} is called the {\em Canonical Polyadic  Decomposition} (CPD) of  $\mathcal T$ and
$R$  is called the rank of the tensor $\mathcal T$ (denoted by $r_{\mathcal T}$).

We write \eqref{eqintro2} as $\mathcal T=[\mathbf A,\mathbf B,\mathbf C]_R$, where
the matrices
$\mathbf A :=\left[\begin{matrix}\mathbf a_1&\dots&\mathbf a_R\end{matrix}\right] \in\mathbb R^{I\times R}$, $\mathbf B :=\left[\begin{matrix}\mathbf b_1&\dots&\mathbf b_R\end{matrix}\right]\in\mathbb R^{J\times R}$ and $\mathbf C :=\left[\begin{matrix}\mathbf c_1&\dots&\mathbf c_R\end{matrix}\right]\in\mathbb R^{K\times R}$
are called the  {\em first}, {\em second} and {\em third factor matrix} of $\mathcal T$, respectively.

Obviously,   $\mathbf a\circ\mathbf b\circ\mathbf c$ has frontal slices $\mathbf a\mathbf b^Tc_1,\dots,\mathbf a\mathbf b^Tc_K\in\mathbb R^{I\times J}$.
Hence, \eqref{eqintro2} is equivalent to the system of matrix identities
\begin{equation}
\mathbf T_k=\sum\limits_{r=1}^R\mathbf a_r\mathbf b_r^Tc_{kr}=\mathbf A\Diag{c^k}\mathbf B^T,\qquad 1\leq k\leq K,\label{eq:matranalog1}
\end{equation}
where  $\mathbf c^k$ denotes the $k$-th column of the matrix $\mathbf C^T$ and $\Diag{c^k}$ denotes a square diagonal matrix with the elements of the vector $\mathbf c^k$ on the main diagonal.

For a  matrix $\mathbf T=[\mathbf t_1\ \cdots\ \mathbf t_J]$, we follow the convention that $\vect{\mathbf T}$ denotes
the column vector obtained by stacking the columns of $\mathbf T$ on top of one another, i.e.,
$
    \textup{vec}(\mathbf T) =
    \left[
    \begin{matrix}
    \mathbf t_1^T& \dots& \mathbf t_J^T
    \end{matrix}
    \right]^T
$. 
The matrix $\Matr{\mathcal T}:=
\left[\begin{matrix}\vect{\mathbf T_1^T}&\dots&\vect{\mathbf T_K^T}\end{matrix}\right]\in\mathbb R^{IJ\times K}
$
is called {\em the matricization  or matrix unfolding} of $\mathcal T$.
The inverse operation is called {\em tensorization}: if  $\mathbf X$ is an $IJ\times K$ matrix, then
$\UnMatr{\mathbf X,I,J}$ is the $I\times J\times K$ tensor such that
 $\Matr{\mathcal T}=\mathbf X$.
From the well-known formula
\begin{equation}
\vect{\mathbf A \Diag{d} \mathbf B^T}=
(\mathbf B\odot \mathbf A)\mathbf d, \qquad \mathbf d\in\mathbb R^R\label{eqmatrtovec}
\end{equation}
it follows that
\begin{equation}
\Matr{\mathcal T}:=
\left[\begin{matrix}(\mathbf A\odot \mathbf B)\mathbf c^1&\dots&(\mathbf A\odot \mathbf B)\mathbf c^K\end{matrix}\right]=
(\mathbf A\odot \mathbf B)\mathbf C^T,\label{eqT_V}
\end{equation}
where ``$\odot$'' denotes the {\em Khatri-Rao product} of matrices:
   $$
   \mathbf A\odot\mathbf B :=
   [\mathbf a_1\otimes\mathbf b_1\ \cdots\ \mathbf a_R\otimes\mathbf b_R]\in \mathbb R^{I J\times R}
   $$
and
``$\otimes$'' denotes the {\em Kronecker product}:
$
\mathbf{a}  \otimes \mathbf{b} =
[a_1b_1\dots a_1b_J\ \dots\  a_Ib_1\dots a_Ib_J]^T.
$

It is clear that in (\ref{eqintro2}) the rank-1 terms can be arbitrarily permuted and that vectors within the same rank-1 term can be arbitrarily scaled provided the overall rank-1 term remains the same. The CPD of a tensor {\em is unique } when it is only subject to these trivial indeterminacies.

\subsection{Problem statement}
The CPD was introduced by F. Hitchcock in \cite{Hitchcock}
  and was later referred to as {\em Canonical Decomposition} (Candecomp) \cite{1970_Carroll_Chang}, {\em Parallel Factor Model} (Parafac) \cite{Harshman1970,1994HarshmanLundy}, and {\em Topographic Components Model}  \cite{1988Topographic}. We refer to the overview papers \cite{KoldaReview,2009Comonetall,Lieven_ISPA,LievenCichocki2013}, the books \cite{Kroonenberg2008,smilde2004multi} and the references therein for background and applications in  Signal Processing, Data Analysis, Chemometrics, and Psychometrics.

Note that in applications one most often deals with  a perturbed version of \eqref{eqintro2}:
  $$
\widehat{\mathcal T}= \mathcal T+\mathcal N=[\mathbf A,\mathbf B,\mathbf C]_R+\mathcal N,
 $$
  where $\mathcal N$ is  an unknown noise tensor and $\widehat{\mathcal T}$ is the given tensor.
The factor matrices of $\mathcal T$ are approximated  by
a solution of the  optimization problem
\begin{equation}
\min \|\widehat{\mathcal T}-[\mathbf A,\mathbf B,\mathbf C]_R\|,\qquad \text{ s.t. }\quad
 \mathbf A\in\mathbb R^{I\times R},\ \mathbf B\in\mathbb R^{J\times R},\ \mathbf C\in\mathbb R^{K\times R},
\label{eq:optproblem}
\end{equation}
 where $\|\cdot\|$ denotes  a suitable  (usually Frobenius) norm \cite{Sorber}.

 In this paper we limit ourselves to  the noiseless case.
 We show that under mild conditions on factor matrices  the CPD is unique and can be found algebraically in the following sense:
 the CPD can be computed by using basic operations on matrices, by computing  compound matrices,
 by taking the orthogonal  complement of a subspace, and by computing generalized  eigenvalue decomposition.
 We make connections with concepts  like permanents, mixed discriminants, and compound matrices,
 which have so far received little attention in applied linear algebra but are of interest.
 Our presentation is in terms of real-valued tensors for notational convenience.
 Complex variants are easily obtained by taking into account complex conjugations.

The heart of the algebraic approach is the following straightforward connection between CPD of a two-slice tensor  and Generalized Eigenvalue Decomposition (GEVD) of a matrix pencil.
 Consider  an  $R\times R\times 2$ tensor $\mathcal T=[\mathbf A,\mathbf B,\mathbf C]_R$, where
$\mathbf A$ and $\mathbf B$ are nonsingular matrices and the matrix $\Diag{d}:=\Diag{c^1}\Diag{c^2}^{-1}$ is defined and has distinct diagonal entries.
From the  equations $\mathbf T_k=\mathbf A\Diag{c^k}\mathbf B^T$, $k=1,2$ it follows easily that
$\mathbf A\Diag{d}\mathbf A^{-1}=\mathbf T_1\mathbf T_2^{-1}$ and
$\mathbf B\Diag{d}\mathbf B^{-1}=(\mathbf T_2^{-1}\mathbf T_1)^T$. Hence, the matrix $\Diag{d}$ can be found (up to permutation of its diagonal entries) from the eigenvalue decomposition of $\mathbf T_1\mathbf T_2^{-1}$ or $(\mathbf T_2^{-1}\mathbf T_1)^T$  and the columns of $\mathbf A$ (resp. $\mathbf B$) are  the eigenvectors of $\mathbf T_1\mathbf T_2^{-1}$ (resp. $(\mathbf T_2^{-1}\mathbf T_1)^T$) corresponding to the $R$ distinct eigenvalues $d_1,\dots,d_R$.
Since the matrices $\mathbf A$ and $\mathbf B$ are nonsingular, the matrix $\mathbf C$ can be easily found from \eqref{eqT_V}. More generally, when $\mathbf A$ and $\mathbf B$ have full column rank and $\mathbf C$ does not have collinear columns, $\mathbf A$ and $\mathbf B$ follow from the GEVD of the matrix pencil $(\mathbf T_1,\mathbf T_2)$.
\subsection{Previous results on uniqueness and algebraic algorithms}

We say that an $I\times R$ matrix {\em has  full column rank } if its column rank is $R$, which implies $I \geq R$.
The  following theorem   generalizes the result discussed at the end of the previous subsection. Several variants of this theorem have appeared in the literature \cite{Harshman1972, TenBerge2009,Lieven_ISPA,Leurgans1993,Sands1980,Sanchez1990}. The proof is essentially obtained by picking two slices (or two mixtures of slices) from $\mathcal T$ and computing their GEVD.
\begin{theorem}\label{theoremoralggeneigdec} Let  $\mathcal T=[\mathbf A,\mathbf B,\mathbf C]_R$ and suppose that
$\mathbf A$ and $\mathbf B$ have full column rank and that  $k_{\mathbf C}\geq 2$.
Then
\begin{itemize}
\item[\textup{(i)}]
 $r_{\mathcal T}=R$ and  the CPD of $\mathcal T$ is unique;
 \item[\textup{(ii)}] the CPD of $\mathcal T$ can be found algebraically.
 \end{itemize}
\end{theorem}
In Theorem \ref{theoremoralggeneigdec} the third  factor matrix plays a different role than the first and the second  factor matrices. Obviously, the theorem still holds when $\mathbf A$, $\mathbf B$, $\mathbf C$ are permuted. In the sequel we  will present only one  version of results. Taking this  into account, we may say that the following result is stronger than
Theorem \ref{theoremoralggeneigdec}.
\begin{theorem}\label{Theorem1.12}
Let  $\mathcal T=[\mathbf A,\mathbf B,\mathbf C]_R$,  $r_{\mathbf C}=R$, and suppose
that $\mathcal C_2(\mathbf A)\odot \mathcal C_2(\mathbf B)$   has full column rank.
Then
\begin{itemize}
\item[\textup{(i)}]
 $r_{\mathcal T}=R$ and  the CPD of $\mathcal T$ is unique  \cite{DeLathauwer2006,JiangSid2004};
 \item[\textup{(ii)}] the CPD of $\mathcal T$ can be found algebraically \cite{DeLathauwer2006}.
 \end{itemize}
\end{theorem}
Computationally, we may obtain from $\mathcal T$ a partially symmetric tensor $\mathcal W$ that has CPD $\mathcal W=[\mathbf C^{-T},\mathbf C^{-T},\mathbf M]_R$ in which both $\mathbf C^{-T}$ and $\mathbf M$ have full column rank and work as in Theorem \ref{theoremoralggeneigdec} to obtain $\mathbf C^{-T}$. The matrices $\mathbf A$ and $\mathbf B$ are subsequently easily obtained from \eqref{eqT_V}.

Also, some algorithms for symmetric CPD have been obtained in the context of algebraic geometry. We refer to \cite{Oeding2013, Landsberg} and references therein. Further, algebraic algorithms have been obtained for CPDs in which factor matrices are subject to constraints (such as orthogonality and Vandermonde) \cite{Mikaelorth, MikaelVDM}.

Our discussion concerns unsymmetric CPD without constraints. Results for the partially and fully symmetric case may be obtained by setting two or all three factor matrices equal to each other, respectively.

In the remaining part of this subsection we present some results on the uniqueness of the CPD.
These results will guarantee CPD uniqueness under the conditions for which we will derive algebraic algorithms.
For more general results on uniqueness we refer to \cite{PartI,PartII}.
The following result was obtained by J. Kruskal, which is little known. We present the compact version from \cite{PartII}.
Corollary \ref{theoremKruskal} presents what is widely known as ``Kruskal's condition'' for CPD uniqueness.
\begin{theorem}\cite[Theorem 4b, p. 123]{Kruskal1977},\cite[Corollary 1.29]{PartII}\label{theoremKruskalnew1}
Let $\mathcal T=[\mathbf A,\mathbf B,\mathbf C]_R$. Suppose that  
 \begin{equation}\label{eqtwomatrK2}
k_{\mathbf A}+r_{\mathbf B}+r_{\mathbf C}                                         \geq 2R+2 \ \ \text{ and } \ \
\min (r_{\mathbf C}+k_{\mathbf B},k_{\mathbf C}+r_{\mathbf B})\geq\ R+2.
 \end{equation}
Then $r_{\mathcal T}=R$ and the CPD of  tensor $\mathcal T$  is unique.
\end{theorem}
\begin{corollary}\cite[Theorem 4a, p. 123]{Kruskal1977}\label{theoremKruskal}
Let $\mathcal T=[\mathbf A,\mathbf B,\mathbf C]_R$  and let
\begin{equation}
k_{\mathbf A}+k_{\mathbf B}+k_{\mathbf C}\geq 2R+2.
 \label{Kruskal}
 \end{equation}
Then $r_{\mathcal T}=R$ and the CPD of  $\mathcal T=[\mathbf A,\mathbf B,\mathbf C]_R$ is  unique.
 \end{corollary}

 In \cite{PartI,PartII} the authors obtained new sufficient conditions expressed in terms of compound matrices.
We will use the following result.
\begin{theorem}\cite[Corollary 1.25]{PartII}\label{Theorem:intro unique}
Let $\mathcal T=[\mathbf A, \mathbf B, \mathbf C]_R$ and
$m:=R-r_{\mathbf C}+2$.
Suppose that
\begin{gather}
\max(\min(k_{\mathbf A},k_{\mathbf B}-1),\ \min(k_{\mathbf A}-1,k_{\mathbf B}))+k_{\mathbf C}\geq R+1,\label{cond(i)}\\
\mathcal C_m(\mathbf A)\odot \mathcal C_m(\mathbf B)   \text{ has full column rank.}\label{cond(ii)}
\end{gather}
Then $r_{\mathcal T}=R$ and the CPD of  tensor $\mathcal T$  is unique.
\end{theorem}

Since the $k$-rank of a matrix cannot exceed its rank (and a fortiori not its number of columns), condition \eqref{Kruskal} immediately implies conditions \eqref{eqtwomatrK2} and  \eqref{cond(i)}.
It was shown in \cite{PartII} that \eqref{eqtwomatrK2} implies  \eqref{cond(ii)}  for $m=R-r_{\mathbf C}+2$.
Thus, Theorem \ref{Theorem:intro unique} guarantees the uniqueness of the CPD under milder conditions than Theorem \ref{theoremKruskalnew1}.
Note also that  statement \textup{(i)} of  Theorem \ref{Theorem1.12} is the special case  of Theorem \ref{Theorem:intro unique} obtained for $r_{\mathbf C}=R$, i.e., when one of the factor matrices has full column rank.
\subsection{New results}\label{subsection1.4}
To simplify the presentation and without loss of generality we will assume throughout the paper
that the third dimension of the tensor $\mathcal T=[\mathbf A,\mathbf B,\mathbf C]_R$ coincides with $r_{\mathbf C}$. (This can  always be achieved in a ``dimensionality reduction'' step: if the  columns of a matrix $\mathbf V$ form an orthonormal basis of the row space of $\Matr{\mathcal T}$ and the matrix $\mathbf A\odot \mathbf B$ has full column rank (as is always the case in the paper), then $r_{\mathbf C}=r_{\Matr{\mathcal T}}=r_{\mathbf V^T\Matr{\mathcal T}}=r_{\mathbf V^T\mathbf C}$,  and by \eqref{eqT_V},
the matrix $\Matr{\mathcal T}\mathbf V=(\mathbf A\odot \mathbf B)\mathbf C^T\mathbf V$ has $r_{\mathbf C}$ columns, which means that the third dimension of the tensor
$\mathcal T_{\mathbf V}:=\UnMatr{\Matr{\mathcal T}\mathbf V ,I,J}$ is equal to $r_{\mathbf C}$; if the CPD $\mathcal T_{\mathbf V}=[\mathbf A,\mathbf B,\mathbf V^T\mathbf C]_R$ has been computed, then
the matrix $\mathbf C$ can be recovered as $\mathbf C=\mathbf V(\mathbf V^T\mathbf C)$).

The following theorems are the main results of the paper. In all cases we will reduce the computation to the situation as in Theorem \ref{theoremoralggeneigdec}.
\begin{theorem}\label{theorem: main}
Let $\mathcal T=[\mathbf A, \mathbf B, \mathbf C]_R$,
$m:=R-r_{\mathbf C}+2$.
Suppose that $k_{\mathbf C}=r_{\mathbf C}$ and  that
\eqref{cond(ii)} holds.
Then
\begin{itemize}
\item[\textup{(i)}]
 $r_{\mathcal T}=R$ and  the CPD of $\mathcal T$ is unique;
 \item[\textup{(ii)}] the CPD of $\mathcal T$ can be found algebraically.
 \end{itemize}
\end{theorem}
Theorem \ref{corrolary: main1} generalizes Theorem \ref{theorem: main} to case where possibly $k_{\mathbf C}<r_{\mathbf C}$. The more general situation for
$\mathbf C$ is accommodated by tightening the condition on $\mathbf A$ and $\mathbf B$. (Indeed, \eqref{cond(secondthmii)} is more restrictive than \eqref{cond(ii)} when
$n>m$.) The proof of Theorem \ref{corrolary: main1} is simple; we essentially consider a $k_{\mathbf C}$-slice subtensor
$\bar{\mathcal T}=[\mathbf A,\mathbf B,\bar{\mathbf C}]_R$ for which $k_{\bar{\mathbf C}}=r_{\bar{\mathbf C}}$, so that Theorem \ref{theorem: main} applies. (Actually,
to guarantee that $k_{\bar{\mathbf C}}=r_{\bar{\mathbf C}}$, we  consider a random slice-mixture.)
\begin{theorem}\label{corrolary: main1}
Let $\mathcal T=[\mathbf A, \mathbf B, \mathbf C]_R$,
$n:=R-k_{\mathbf C}+2$.
Suppose that
\begin{equation}
\mathcal C_n(\mathbf A)\odot \mathcal C_n(\mathbf B)   \text{ has full column rank.}\label{cond(secondthmii)}
\end{equation}
Then
\begin{itemize}
\item[\textup{(i)}]
 $r_{\mathcal T}=R$ and  the CPD of $\mathcal T$ is unique;
 \item[\textup{(ii)}] the CPD of $\mathcal T$ can be found algebraically.
 \end{itemize}
\end{theorem}
We also obtain the following corollaries.
\begin{corollary}\label{likegenKruskal}
Let $\mathcal T=[\mathbf A,\mathbf B,\mathbf C]_R$. Suppose that
 \begin{equation}\label{eqtwomatrK2kc=rc}
k_{\mathbf A}+r_{\mathbf B}+k_{\mathbf C}                                         \geq 2R+2,\ \ \text{and}\ \
k_{\mathbf B}+k_{\mathbf C}    \geq\ R+2.
 \end{equation}
Then $r_{\mathcal T}=R$ and the CPD of  tensor $\mathcal T$  is unique and can be found algebraically.
\end{corollary}
\begin{corollary}\label{corrolary: main2}
Let $\mathcal T=[\mathbf A,\mathbf B,\mathbf C]_R$  and let $k_{\mathbf A}+k_{\mathbf B}+k_{\mathbf C}\geq 2R+2$.
Then  the CPD of $\mathcal T$ is unique and  can be found algebraically.
\end{corollary}

Let us further explain how the theorems that we have formulated so far, relate to one another. First,  we  obviously
have that $n=R-k_{\mathbf C}+2\geq R-r_{\mathbf C}+2=m$. Next, the following implications were proved in  \cite{PartI}:
\begin{equation}\label{scheme}
\begin{tikzcd}
\eqref{eqtwomatrK2kc=rc} \arrow[Rightarrow]{r}\arrow[Rightarrow]{d}{\text{trivial}} & \eqref{cond(secondthmii)} \arrow[Rightarrow]{r}\arrow[Rightarrow]{d}& \min(k_{\mathbf A},k_{\mathbf B})\geq n
\arrow[Rightarrow,
start anchor={[xshift=-3ex]},
end anchor={[xshift=-3ex]}]{d}
\arrow[Rightarrow]{r}{\text{trivial}}& \eqref{cond(i)}  \\
 \eqref{eqtwomatrK2} \arrow[Rightarrow]{r}                        & \eqref{cond(ii)} \arrow[Rightarrow]{r}                     & \min(k_{\mathbf A},k_{\mathbf B})\geq m \arrow[Rightarrow]{u}[right]{\text{if }k_{\mathbf C}=r_{\mathbf C}\text{ (trivial)}}
\end{tikzcd}
\end{equation}
The first thing that follows from scheme \eqref{scheme} is that Theorem \ref{corrolary: main1}  is indeed more general than Corollary \ref{likegenKruskal}.
Corollary \ref{corrolary: main2} follows trivially from Corollary \ref{likegenKruskal}.
Next, it appears that the conditions of Theorems
 \ref{theorem: main}--\ref{corrolary: main1} are more restrictive than the conditions of Theorem  \ref{Theorem:intro unique}. Also, the conditions of Corollary
 \ref{likegenKruskal} are more restrictive than the conditions of Theorem \ref{theoremKruskalnew1}.
Hence,  we immediately obtain the uniqueness of the CPD in
Theorems   \ref{theorem: main}--\ref{corrolary: main1} and Corollary \ref{likegenKruskal}.
Consequently,  we can limit ourselves to the derivation of the algebraic algorithms.
\subsection{Organization}We  now explain how the paper is organized. Let $\mathcal T=[\mathbf A,\mathbf B,\mathbf C]_R\in\mathbb R^{I
\times J\times K}$ with $k_{\mathbf C}=K$, implying $K\leq R$.
In  the first phase of our algorithms,   we find up to column permutation and scaling the $K\times C^{K-1}_R$ matrix
$\mathcal B(\mathbf C)$ defined by
\begin{equation}
\mathcal B(\mathbf C):=\mathbf L\mathcal C_{K-1}(\mathbf C),\label{def:B(C)}
\end{equation}
where
\begin{equation}
\mathbf L:=\left[\begin{array}{rrrc} 0     & 0     & \dots               &(-1)^{K-1}\\
                                     \vdots& \vdots& \iddots           &\vdots\\
                                0     & -1    & \dots               &0\\
                                1     & 0     & \dots                  &0
                 \end{array}\right].
                 \label{eq:matrixL}
\end{equation}
The matrix $\mathcal B(\mathbf C)$ can be considered as an unconventional variant of the inverse of $\mathbf C$:
\begin{align}
&\text{every column of } \mathcal B(\mathbf C) \text{ is orthogonal to exactly } K-1 \text{ columns of }\mathbf C, \tag{P1}\label{P1}\\
&\begin{multlined}
 \text{any vector that is orthogonal to exactly } K-1 \text{ columns of } \mathbf C \\
 \qquad\qquad\qquad\text{ is proportional to a column of } \mathcal B(\mathbf C),
\end{multlined} \tag{P2}\label{P2}
\\
&\text{ every column of } \mathbf C \text{ is orthogonal to exactly } C^{K-2}_{R-1} \text{ columns of } \mathcal B(\mathbf C),\tag{P3}\label{P3}\\
&\begin{multlined}
\text{any vector that is orthogonal to exactly } C^{K-2}_{R-1} \text{ columns of } \mathcal B(\mathbf C)\\
\text{ is proportional to a column of } \mathbf C.
 \end{multlined}\tag{P4}\label{P4}
\end{align}
Recall that every column of the classical Moore-Penrose pseudo-inverse $\mathbf C^\dagger\in\mathbb R^{R\times K}$
is orthogonal to exactly $K-1$ {\em rows} of $\mathbf C$ and vice-versa.
The equality $\mathbf C\mathbf C^\dagger=\mathbf I_K$ works along the ``long'' dimension of $\mathbf C$.
If $\mathbf C^\dagger$ is known, then $\mathbf C$ may easily be found by pseudo-inverting again, $\mathbf C=(\mathbf C^\dagger)^\dagger$.
The interaction with $\mathcal B(\mathbf C)$ takes place along the ``short'' dimension of $\mathbf C$, and this complicates things.
Nevertheless, it is also possible to reconstruct $\mathbf C$
from $\mathcal B(\mathbf C)$. 
In the second and third phase of our algorithms we use $\mathcal B(\mathbf C)$ to compute  CPD.
The following two properties of $\mathcal B(\mathbf C)$ will be crucial for our derivation.
\begin{proposition}\label{prop:newB(C)}
Let $\mathbf C\in\mathbb R^{K\times R}$ and $k_{\mathbf C}=K$.
Then
\begin{itemize}
\item[\textup{(i)}]
$\mathcal B(\mathbf C)$ has no proportional columns, that is $k_{\mathcal B(\mathbf C)}\geq 2$.
\item[\textup{(ii)}] the matrices
$$
\mathcal B(\mathbf C)^{(m-1)}=
\underbrace{\mathcal B(\mathbf C)\odot\dots\odot \mathcal B(\mathbf C)}_{m-1},\qquad
\mathcal B(\mathbf C)^{(m)}=
\underbrace{\mathcal B(\mathbf C)\odot\dots\odot \mathcal B(\mathbf C)}_{m}
$$
 have full column rank for  $m:=R-K+2$.
 \end{itemize}
\end{proposition}

Sections  \ref{sec:compandpermcomp}--\ref{sec:multilinear maps} contain auxiliary results of which several are interesting in their own right.
In Subsection \ref{subsection 2.1} we recall the properties of compound matrices, provide an intuitive understanding of
 properties \eqref{P1}--\eqref{P4} and Propositions  \ref{prop:newB(C)}, and discuss the reconstruction of 
 $\mathbf C$ from  $\mathcal B(\mathbf C)$.
 (Since the proofs of properties \eqref{P1}-\eqref{P4}
 and Proposition \ref{prop:newB(C)} are rather long and technical, they are included in the supplementary materials.)
  In Subsections \ref{subsection 2.4}--\ref{subsection 2.5} we study
variants of permanental compound matrices.
 Let the columns of the $K^m$-by-$C^m_R$ matrix $\mathcal R_m(\mathbf C)$
be  equal to  the vectorized symmetric parts of the tensors
 $\mathbf c_{i_1}\circ\dots\circ\mathbf c_{i_m}$, $1\leq i_1<\dots<i_m\leq R$ and let 
 $\textup{range}(\pi_S)$ denote a subspace of $\mathbb R^{K^m}$ that consists of 
 vectorized versions of $m$-th order $K\times\dots\times K$ symmetric tensors, yielding $\dim \textup{range}(\pi_S)=C^m_{K+m-1}$.
 We prove
 \begin{equation}
\text{Proposition } \ref{proposition:2.23}\ \textup{(iii)}:\
\ker\left(\mathcal{R}_m(\mathbf C)^T\upharpoonright_{\textup{range}(\pi_S)}\right)= \textup{range} (\mathcal B(\mathbf C)^{(m)}),\label{eq:mostimportnat}
\end{equation}
where the notation $\mathcal{R}_m(\mathbf C)^T\upharpoonright_{\textup{range}(\pi_S)}$
means that we let the matrix $\mathcal{R}_m(\mathbf C)^T$ act only 
on vectors from $\textup{range}(\pi_S)$, i.e., on $K^m\times 1$ vectorized versions of $K\times\dots\times K$ symmetric tensors. Computationally,   the subspace $\ker\left(\mathcal{R}_m(\mathbf C)^T\upharpoonright_{\textup{range}(\pi_S)}\right)$ is the intersection of the subspaces $\ker(\mathcal{R}_m(\mathbf C)^T)$ and $\textup{range}(\pi_S)$.

In \S \ref{sec:multilinear maps} we introduce polarized compound matrices --- a notion closely related to
the rank detection mappings in \cite{DeLathauwer2006, DimLievenLL1}.
 The entries of polarized compound matrices are mixed  discriminants \cite{determinantfunctions, Bapat1989, Aleksandrov1938}.
Using polarized compound matrices we construct a
$C^m_IC^m_J\times K^m$ matrix $\RTm{m}$ from the given tensor $\mathcal T$ such that
\begin{equation}
\RTm{m}=\left[\mathcal C_m(\mathbf A)\odot\mathcal C_m(\mathbf B)\right]\mathcal R_m(\mathbf C)^T.
\label{eq:mainformula}
\end{equation}
Assuming that $\mathcal C_m(\mathbf A)\odot\mathcal C_m(\mathbf B)$ has full column rank and combining \eqref{eq:mostimportnat} with
\eqref{eq:mainformula} we find  the space generated by the columns of the matrix $\mathcal B(\mathbf C)^{(m)}$:
\begin{equation}
\ker\left(\RTm{m}\upharpoonright_{\textup{range}(\pi_S)}\right)=
\ker\left(\mathcal{R}_m(\mathbf C)^T\upharpoonright_{\textup{range}(\pi_S)}\right)=
\textup{range} (\mathcal B(\mathbf C)^{(m)}).
\label{eq:intro1}
\end{equation}

In \S  \ref{Section4} we combine all results to obtain Theorems \ref{theorem: main}--\ref{corrolary: main1} and we present two algebraic CPD  algorithms. 
Both new algorithms contain the same first phase in which we
find a matrix $\mathbf F$ that coincides with $\mathcal B(\mathbf C)$ up to column permutation and scaling. 
This first phase of the algorithms relies on  key formula
\eqref{eq:intro1}, which makes a link between the known matrix
$\RTm{m}$, constructed from $\mathcal T$, and the unknown matrix
$\mathcal B(\mathbf C)$. We work as follows. We construct the matrix
$\RTm{m}$  and compute the vectorized symmetric tensors in its kernel. We stack a basis of 
$\ker\left(\RTm{m}\upharpoonright_{\textup{range}(\pi_S)}\right)$ as columns of a matrix
$\Matr{\mathcal W}\in\mathbb R^{K^m\times C^{K-1}_R}$, with which we associate a
$K\times K^{m-1}\times C^{K-1}_R$ tensor $\mathcal W$. 
From Proposition \ref{prop:newB(C)} and Theorem \ref{theoremoralggeneigdec} it follows that the  CPD $\mathcal W=[\mathcal  B(\mathbf C)$, $\mathcal  B(\mathbf C)^{(m-1)}, \mathbf M]_{C^{K-1}_R}$ can be found algebraically.  This allows us to find a matrix $\mathbf F$ that coincides with $\mathcal B(\mathbf C)$ up to column permutation and scaling. In the second and third phase of the  first algorithm we  find the matrix $\mathbf C$  and the matrices $\mathbf A$ and $\mathbf B$, respectively.
For finding $\mathbf C$, we resort to properties \eqref{P3}--\eqref{P4}. Full exploitation of the structure has combinatorial complexity and is
infeasible unless the dimensions of the tensor are relatively small. As an alternative, in the second algorithm we first find
 the matrices $\mathbf A$ and $\mathbf B$ and then we find the matrix $\mathbf C$.
This is done as follows. We construct  the new $I\times J\times C^{K-1}_R$ tensor $\mathcal V$ with the matrix unfolding
$\Matr{\mathcal V}:=\Matr{\mathcal T}\mathbf F=
(\mathbf A\odot \mathbf B)\mathbf C^T\mathbf F$.
We find subtensors of $\mathcal V$ such that
each subtensor has dimensions $I\times J\times 2$ and its CPD can be found algebraically. Full exploitation of the structure yields
$C^m_RC^2_m$ subtensors. From the CPD of the subtensors  we simultaneously obtain the columns of $\mathbf A$ and $\mathbf B$, and finally we set  $\mathbf C=\left((\mathbf A\odot\mathbf B)^\dagger\Matr{\mathcal T}\right)^T$.

We conclude the paper with two examples.
In the first example we demonstrate how the  algorithms work for
a  $4\times 4\times 4$ tensor of rank $5$ for which  $k_{\mathbf A}=k_{\mathbf B}=3$.
In the second  example we consider
a generic $6\times 6\times 7$ tensor of rank $9$ and compare
the complexity of algorithms.
Note that in  neither case the uniqueness of the CPDs follows from Kruskal's Theorem \ref{theoremKruskalnew1}.
\subsection{Link with \cite{DeLathauwer2006}}
Our overall derivation generalizes ideas from \cite{DeLathauwer2006} ($K=R$).
To conclude the introduction, we recall the
CPD  algorithm   from \cite{DeLathauwer2006} using our notations. We have
  $K=R$, which implies  $m=2$.
  First, we construct the  $C^2_IC^2_J\times R^2$ matrix $\mathbf R_2(\mathcal T)$
   whose $((i-1)R+j)$-th column is computed as
  \begin{equation*}
   \text{Vec}\left(\ \mathcal C_2(\mathbf T_i+\mathbf T_j) - \mathcal C_2(\mathbf T_i) -\mathcal C_2(\mathbf T_j)\ \right),\quad 1\leq i\leq j\leq R, 
   \end{equation*}
   where  $\mathbf T_1,\dots,\mathbf T_R\in\mathbb R^{I\times J}$ denote the frontal slices of $\mathcal T$. The  entries of the $((i-1)R+j)$-th column of $\mathbf R_2(\mathcal T)$  can be identified with  the  $C^2_IC^2_J$ nonzero entries of the   $I\times I\times J\times J$ tensor $\mathcal P_{ij}$ \cite[p. 648]{DeLathauwer2006}.
     Then we find a basis $\mathbf w_1,\dots,\mathbf w_R\in\mathbb R^{R^2}$ of 
  $E:=\ker\left(\RTm{2}\upharpoonright_{\textup{range}(\pi_S)}\right)$ and set 
  $\mathbf W=[\mathbf w_1\ \dots\ \mathbf w_R]$. 
  We note that $E$ can be computed as  the intersection of the subspaces  $\ker(\RTm{2})$ and $\textup{range}(\pi_S)$, where
  $\textup{range}(\pi_S)$  consists of  vectorized  versions of symmetric $R\times R$ matrices.
  In \cite{DeLathauwer2006}, the subspace $E$ is generated by the vectors in $\textup{range}(\pi_S)$
  that yield a zero linear combination of the $R^2$  tensors $\mathcal P_{ij}$. 
           In the next step we recover (up to column permutation and scaling) $\mathbf C$ from  $E$. This is done as follows.
   By \eqref{P3}--\eqref{P4}, the columns of $\mathcal B(\mathbf C)$ are proportional to the columns of $\mathbf C^{-T}$,
   i.e., $\mathcal B(\mathbf C)^T$ is equal to the inverse of $\mathbf C$ up to column permutation and scaling.
Hence, by \eqref{eq:intro1}, $\textup{range}(\mathbf W)=\textup{range}(\mathbf C^{-T}\odot \mathbf C^{-T})$.
Hence, there exists a nonsingular matrix $\mathbf M$ such that  $\mathbf W=
\left(\mathbf C^{-T}\odot \mathbf C^{-T}\right)\mathbf M^T$.
Therefore, by \eqref{eqT_V},
$\mathcal W=[\mathbf C^{-T},\mathbf C^{-T},\mathbf M]_R$, where
$\mathcal W$ denotes the $R\times R\times R$ tensor such that $\mathbf W=\Matr{\mathcal W}$. Since all factor matrices of $\mathcal W$ have full column rank, the CPD of $\mathcal W$ can be computed algebraically.
Thus, we can find $\mathbf C^{-T}$ (and hence, $\mathbf C$) up to column permutation and scaling. Finally, the matrices $\mathbf A$ and $\mathbf B$ can now be easily found from $\Matr{\mathcal T}\mathbf C^{-T}=
\mathbf A\odot \mathbf B$ using the fact that the columns of
$\mathbf A\odot \mathbf B$ are vectorized rank-$1$ matrices.

\section{Matrices formed by determinants and permanents of submatrices of a given matrix}\label{sec:compandpermcomp}
 Throughout the paper we will use the following multi-index notations.
Let $i_1,\dots,i_k$ be integers. Then $\{i_1,\dots,i_k\}$  denotes the set with elements $i_1,\dots,i_k$ (the order does not matter) and $(i_1,\dots,i_k)$ denotes a $k$-tuple (the order is important).
 Let
 \begin{eqnarray*}
S^k_n&=&\{(i_1,\dots,i_k):\ 1\leq i_1<i_2<\dots<i_k\leq n\},\\
Q^k_n&=&\{(i_1,\dots,i_k):\ 1\leq i_1\leq i_2\leq \dots\leq i_k\leq n\},\\
R^k_n&=&\{(i_1,\dots,i_k):\ i_1,\dots,i_k\in\{1,\dots,n\}\}.
\end{eqnarray*}
It is well known that
$
\textup{card}\ S^k_n=C^k_n$ ,
$\textup{card}\ Q^k_n=C^k_{n+k-1}$, and
$\textup{card}\ R^k_n=n^k$.
We assume that the elements  of  $S^k_n$, $Q^k_n$, and $R^k_n$ are ordered lexicographically.
In the sequel we will both use indices taking values in $\{1,2, \dots, C^{k}_n\}$ (resp.  $\{1,2, \dots, C^{k}_{n+k-1}\}$ or $\{1,2, \dots, n^k\}$)
and multi-indices taking values in $S_{n}^k$ (resp. $Q_{n}^k$ or $R_{n}^k$).
For example,
\begin{gather*}
S_2^2=\{(1,2)\},\quad  Q_2^2=\{(1,1),(1,2),(2,2)\},\quad R_2^2=\{(1,1),(1,2),(2,1),(2,2)\},\\
S_2^2(1)=Q_2^2(2)=R_2^2(2),\quad Q_2^2(3)=R_2^2(4).
\end{gather*}
Let also $P_{\{j_1,\dots,j_n\}}$ denote the set of all permutations of the set $\{j_1,\dots,j_n\}$.
We follow the convention that if some of  $j_1,\dots,j_n$ coincide, then
the set $P_{\{j_1,\dots,j_n\}}$ contains identical elements, yielding $\textup{card}\ P_{\{j_1,\dots,j_n\}}=n!$.
For example, $
P_{\{1,2,2\}}=\{\{1,2,2\},\{1,2,2\},\{2,1,2\},\{2,2,1\},\{2,1,2\},\{2,2,1\}\}.
$
We set $P_n:=P_{\{1,\dots,n\}}$.

Let $\mathbf A\in \mathbb R^{m\times n}$. Throughout the paper
$\mathbf A ((i_1,\dots,i_k), (j_1,\dots,j_k))$ denotes the submatrix of $\mathbf A$
at the intersection of the $k$ rows with row numbers $i_1,\dots,i_k$ and the $k$ columns with column numbers $j_1,\dots,j_k$.

\subsection{Matrices whose entries are determinants}\label{subsection 2.1}
In this subsection we briefly discuss  compound matrices.
The $k$-th compound matrix of a given matrix is formed by $k\times k$ minors of that matrix. We have the following formal definition.
\begin{definition}\label{def:compound matrices} \cite{HornJohnson}
Let $\mathbf A\in \mathbb R^{m\times n}$ and $k\leq\min(m,n)$.
The $C_m^k$-by-$C_n^k$ matrix whose $(i, j)$-th entry is $\det\mathbf A (S_m^k(i), S_n^k(j))$ is
called the $k$-th compound matrix of $\mathbf A$ and is denoted by $\mathcal C_k(\mathbf A)$.
\end{definition}
\begin{example}\label{Example2.2}
Let $\mathbf A=[\mathbf I_3\ \mathbf a]$, where $\mathbf a=[a_1\ a_2\ a_3]^T$.
Then
\begin{align*}
\mathcal C_2(\mathbf A)
=&\bordermatrix[{[]}]{%
&(1,2)&(1,3)&(1,4)&(2,3)&(2,4)&(3,4)\cr
(1,2)&
\Big|
\begin{matrix}
    1& 0\\
    0& 1
    \end{matrix}
\Big|
&
\Big|
\begin{matrix}
    1& 0\\
    0& 0
    \end{matrix}
\Big|
&
\Big|
\begin{matrix}
   1 & a_1\\
   0 & a_2
    \end{matrix}
\Big|
&
\Big|
\begin{matrix}
    0& 0\\
    1& 0
    \end{matrix}
\Big|
&
\Big|
\begin{matrix}
    0& a_1\\
    1& a_2
    \end{matrix}
\Big|
&
\Big|
\begin{matrix}
    0& a_1\\
    0& a_2
    \end{matrix}
\Big|
\cr
(1,3)&
\Big|
\begin{matrix}
    1& 0\\
    0& 0
    \end{matrix}
\Big|
&
\Big|
\begin{matrix}
    1& 0\\
    0& 1
    \end{matrix}
\Big|
&
\Big|
\begin{matrix}
    1& a_1\\
    0& a_3
    \end{matrix}
\Big|
&
\Big|
\begin{matrix}
    0& 0\\
    0& 1
    \end{matrix}
\Big|
&
\Big|
\begin{matrix}
    0& a_1\\
    0& a_3
    \end{matrix}
\Big|
&
\Big|
\begin{matrix}
    0& a_1\\
    1& a_3
    \end{matrix}
\Big|
\cr
(2,3)&
\Big|
\begin{matrix}
    0& 1\\
    0& 0
    \end{matrix}
\Big|
&
\Big|
\begin{matrix}
    0& 0\\
    0& 1
    \end{matrix}
\Big|
&
\Big|
\begin{matrix}
    0& a_2\\
    0& a_3
    \end{matrix}
\Big|
&
\Big|
\begin{matrix}
    1& 0\\
    0& 1
    \end{matrix}
\Big|
&
\Big|
\begin{matrix}
    1& a_2\\
    0& a_3
    \end{matrix}
\Big|
&
\Big|
\begin{matrix}
    0& a_2\\
    1& a_3
    \end{matrix}
\Big|
}\\
=& \left[
\begin{array}{rrrrrr}
1& 0&a_2&0&-a_1&0\\
0&1&a_3&0&0&-a_1\\
0&0&0&1&a_3&-a_2
\end{array}
\right].
\end{align*}
\end{example}
Definition \ref{def:compound matrices}  immediately implies the following lemma.
\begin{lemma} \label{compoundprop1}
Let $\mathbf A\in\mathbb R^{I\times R}$ and $k\leq \min(I,R)$. Then
\begin{itemize}
\item[\textup{(1)}]
 $\mathcal C_k(\mathbf A)$ has one or more zero columns if and only if
$k > k_{\mathbf A}$;
\item[\textup{(2)}]
$\mathcal C_k(\mathbf A)$ is equal to the zero matrix if and only if
$k > r_{\mathbf A}$;
\item[\textup{(3)}]
$\mathcal C_k(\mathbf A^T) = (\mathcal C_k(\mathbf A))^T$.
\end{itemize}
\end{lemma}
PD representation \eqref{eq:matranalog1} will make us need compound matrices of diagonal matrices.
\begin{lemma}\label{Lemma2.3}
Let  $\mathbf d\in\mathbb R^R$, let $\omega(\mathbf d)$ denote the number of nonzero entries of  $\mathbf d$,  $k\leq R$, and let
$
\hatdSmR{k}{R}:=
[d_1\cdots d_k\ \
d_1\cdots d_{k-1}d_{k+1}\ \
\dots$\ \ $d_{R-k+1}\cdots d_R]^T\in\mathbb R^{C^k_R}.
$
Then
\begin{itemize}
\item[\textup{(1)}]
$\hatdSmR{k}{R}=\vzero$ if and only if $\omega(\mathbf d)\leq k-1$;
\item[\textup{(2)}]
$\hatdSmR{k}{R}$ has exactly one nonzero entry if and only if $\omega(\mathbf d)=k$;
\item[\textup{(3)}]
$\mathcal C_k(\textup{\text{Diag}}(\mathbf d))=\textup{\text{Diag}}(\hatdSmR{k}{R})$.
\end{itemize}
\end{lemma}
The following result is known as Binet-Cauchy formula.
\begin{lemma}\label{LemmaCompound}\cite[p. 19--22]{HornJohnson}
Let $k$  be a positive integer and let $\mathbf A$ and $\mathbf B$ be
matrices   such that  $\mathcal C_k(\mathbf A)$ and  $\mathcal C_k(\mathbf B)$, are  defined.
Then $\mathcal C_k(\mathbf A\mathbf B^T) =\mathcal  C_k(\mathbf A)\mathcal C_k(\mathbf B^T)$. If additionally $\mathbf d$ is a vector such that
$\mathbf A\textup{\text{Diag}}(\mathbf d)\mathbf B^T$ is defined, then
$
\mathcal C_k(\mathbf A \textup{\text{Diag}}(\mathbf d) \mathbf B^T)=
\mathcal C_k(\mathbf A) \textup{\text{Diag}}(\hatdSmR{k}{R})\mathcal C_k( \mathbf B)^T.
$
\end{lemma}

The goal of the remaining part of this subsection is to provide an intuitive understanding of
properties \eqref{P1}--\eqref{P4} and Proposition \ref{prop:newB(C)}.

Let $K\geq 2$, and let $\mathbf C$ be a $K\times K$ nonsingular matrix.
By Cramer's rule and \eqref{def:B(C)},
the  matrices $\det(\mathbf C)\mathbf C^{-1}$ and  $\mathcal B(\mathbf C)$  are formed by  $(K-1)\times (K-1)$ minors (also known as
cofactors) of $\mathbf C$. It is easy to show that  $\mathcal B(\mathbf C)=(\det(\mathbf C)\mathbf C^{-1})^T\mathbf L$, where
$\mathbf L$ is given by \eqref{eq:matrixL}.
It now trivially follows that every column of $\mathcal B(\mathbf C)$ is a nonzero vector orthogonal to exactly $K-1$ columns of $\mathbf C$. Indeed,
 $$
 \mathbf C^T\mathcal B(\mathbf C)=\mathbf C^T\det(\mathbf C)\mathbf C^{-T}\mathbf L=
 \det(\mathbf C)\mathbf L,
 $$
 which has precisely one non-zero entry in every column.
The inverse statement holds also. Namely, if $\mathbf x$ is a nonzero vector 
that is orthogonal to exactly $K$ ($=C^{K-2}_{K-1}$) columns of  $\mathcal B(\mathbf C)$ (i.e.  $\omega(\mathbf x^T\mathcal B(\mathbf C))\leq 1$), then
$\mathbf x$ is proportional to a column of $\mathbf C$. Indeed,
\begin{equation}\label{eq:condKK}
\begin{split}
\omega(\mathbf x^T\mathcal B(\mathbf C))&=\omega(\mathbf x^T\det(\mathbf C)\mathbf C^{-T}\mathbf L)=
\omega(\mathbf x^T\mathbf C^{-T})=\omega(\mathbf C^{-1}\mathbf x)\leq 1\Leftrightarrow\\
&\mathbf x \ \text{ is proportional to a column of }\ \mathbf C.
\end{split}
\end{equation}
 Properties \eqref{P3}--\eqref{P4}
  generalize \eqref{eq:condKK} for rectangular matrices and imply that, if we know $\mathcal B(\mathbf C)$ up to column permutation and scaling,
then we know $\mathbf C$ up to column  permutation and scaling. This result will be directly used in Algorithm \ref{alg:Main} further: we will
first estimate $\mathcal B(\mathbf C)$ up to column permutation and scaling and then obtain $\mathbf C$ up to column permutation and scaling.
Statements \eqref{P1}--\eqref{P3} are easy to show. Statement \eqref{P4} is more difficult. Since the proofs
are technical, they are given in the supplementary materials.

Let us illustrate properties \eqref{P1}--\eqref{P4} and Proposition \ref{prop:newB(C)} for a rectangular matrix $\mathbf C$ ($K<R$).

\begin{example}\label{example:2.8}
Let
$$
\mathbf C=
\left[
\begin{matrix}
1&0&0&1\\
0&1&0&1\\
0&0&1&1
\end{matrix}
\right],\qquad  \mathbf L=\left[
\begin{array}{rrr}
0&0&1\\
0&-1&0\\
1&0&0
\end{array}
\right],
$$
implying $k_{\mathbf C}=K=3$ and $R=4$. From \eqref{def:B(C)} and Example \ref{Example2.2} it follows that
$$
\mathcal B(\mathbf C)=\mathbf L\mathcal C_2(\mathbf C)=\left[
\begin{array}{rrrrrr}
0&0&0&1&1&-1\\
0&-1&-1&0&0&1\\
1& 0&1&0&-1&0
\end{array}
\right].
$$
One can easily check the statements of properties \eqref{P1}--\eqref{P4} and Proposition \ref{prop:newB(C)}.
Note in particular that exactly $4$ sets of $3$ columns of $\mathcal B(\mathbf C)$  are linearly dependent. The vectors
that are orthogonal  to these sets are proportional to the columns of $\mathbf C$.
\end{example}

In our overall CPD algorithms we will find a matrix $\mathbf F\in\mathbb R^{K\times C^{K-1}_R}$  that
coincides with $\mathcal B(\mathbf C)$ up to column permutation and  scaling.
Properties \eqref{P3}--\eqref{P4} imply the following combinatorial procedure to find  the third factor matrix of $\mathcal T$.
Since the permutation indeterminacy  makes that we do not know beforehand which columns of $\mathbf F$
are orthogonal to which columns of $\mathbf C$,  we need to look for subsets of $C^{K-2}_{R-1}$ columns of $\mathbf F$ that are linearly dependent.
By properties \eqref{P3}--\eqref{P4}, there exist exactly $R$ such subsets. For each subset, the orthogonal complement yields, up to scaling,
a column of $\mathbf C$.

\subsection{Matrices whose entries are permanents}\label{subsection 2.4}\hfill

\begin{definition}\label{def:permanent}
Let $\mathbf A=\left[\begin{matrix}\mathbf a_1&\dots&\mathbf a_n\end{matrix}\right]\in \mathbb R^{n\times n}$. Then the permanent of $\mathbf A$ is defined as
$$
\textup{perm\ }\mathbf A=\overset{+}{|}\mathbf A\overset{+}{|}=\sum\limits_{(l_1,\dots,l_n)\in P_n}a_{1l_1}a_{2l_2}\cdots a_{nl_n}
=\sum\limits_{(l_1,\dots,l_n)\in P_n}a_{l_11}a_{l_22}\cdots a_{l_nn}.
$$
\end{definition}
The definition of the permanent of $\mathbf A$ differs from that of the determinant of $\mathbf A$ in that the signatures of the permutations are not taken into account.
This makes the permanent invariant for column permutations of $\mathbf A$.
The notations $\textup{perm\ }\mathbf A$ and  $\overset{+}{|}\mathbf A\overset{+}{|}$ are due to Minc \cite{Minc1978} and  Muir \cite{Muir1882}, respectively.

We have the following permanental variant of compound matrix.
\begin{definition}\label{def:2.9}\cite{MarcusMinc1961}
Let $\mathbf C\in \mathbb R^{K\times R}$.
The $C_K^m$-by-$C_R^m$ matrix whose $(i, j)$-th entry is $\textup{perm}\ \mathbf C (S_K^m(i), S_R^m(j))$ is
called the $m$-th permanental compound matrix of $\mathbf C$ and is denoted by $\mathcal {PC}_m(\mathbf C)$.
\end{definition}

In our derivation we will also use the following two types of matrices. As far as we know, these do not have a special name.
\begin{definition}\label{def:permanental compound matrices}
Let $\mathbf C\in \mathbb R^{K\times R}$.
The $C_{K+m-1}^m$-by-$C_{R}^m$ matrix whose $(i, j)$-th entry is $\textup{perm\ }\mathbf C (Q_K^m(i), S_R^m(j))$
is denoted by $\mathcal{Q}_m(\mathbf C)$.
\end{definition}
\begin{definition}\label{def:extended permanental compound matrices}
Let $\mathbf C\in \mathbb R^{K\times R}$.
The $K^m$-by-$C_R^m$ matrix whose $(i, j)$-th entry is $\textup{perm\ }\mathbf C (R_K^m(i), S_R^m(j))$
is denoted by $\mathcal{R}_m(\mathbf C)$.
\end{definition}

Note that $\mathcal{Q}_m(\mathbf C)$ is a submatrix of $\mathcal{R}_m(\mathbf C)$, in which the doubles of rows that are due to the permanental
invariance for column permutations, have been removed.

The following lemma makes the connection between
$\mathcal{Q}_m(\mathbf C)^T$  and $\mathcal{R}_m(\mathbf C)^T$ and permanental compound matrices.
\begin{lemma}\label{lemma:permanental compound_equiv2}
Let $\mathbf C=
[\begin{matrix} \mathbf c^1& \dots& \mathbf c^K \end{matrix}]^T\in \mathbb R^{K\times R}$. Then
 $\mathcal{Q}_m(\mathbf C)^T$ (resp. $\mathcal{R}_m(\mathbf C)^T$) has columns $\mathcal{PC}_m([\begin{matrix} \mathbf c^{j_1}& \dots& \mathbf c^{j_m} \end{matrix}])$, where
$(j_1,\dots,j_m)\in Q^m_K$ (resp. $R^m_K$).
\end{lemma}
\begin{example}\label{example:EPC2}
Let $\mathbf C=\left[\begin{matrix}1&2&3\\ 4&5&6\end{matrix}\right]$. Then
\begin{equation*}
\mathcal{R}_2(\mathbf C)=
\bordermatrix[{[]}]{%
&(1,2)&(1,3)&(2,3)\cr
(1,1)&
\overset{+}{\Big|}
\begin{matrix}
    1& 2\\
    1& 2
    \end{matrix}
\overset{+}{\Big|}
&
\overset{+}{\Big|}
\begin{matrix}
    1& 3\\
    1& 3
    \end{matrix}
\overset{+}{\Big|}
&
\overset{+}{\Big|}
\begin{matrix}
    2& 3\\
    2& 3
    \end{matrix}
\overset{+}{\Big|}
\cr
(1,2)&
\overset{+}{\Big|}
\begin{matrix}
    1& 2\\
    4& 5
    \end{matrix}
\overset{+}{\Big|}
&
\overset{+}{\Big|}
\begin{matrix}
    1& 3\\
    4& 6
    \end{matrix}
\overset{+}{\Big|}
&
\overset{+}{\Big|}
\begin{matrix}
    2& 3\\
    5& 6
    \end{matrix}
\overset{+}{\Big|}
\cr
(2,1)&
\overset{+}{\Big|}
\begin{matrix}
    4& 5\\
    1& 2
    \end{matrix}
\overset{+}{\Big|}
&
\overset{+}{\Big|}
\begin{matrix}
    4& 6\\
    1& 3
    \end{matrix}
\overset{+}{\Big|}
&
\overset{+}{\Big|}
\begin{matrix}
    5& 6\\
    2& 3
    \end{matrix}
\overset{+}{\Big|}
\cr
(2,2)&
\overset{+}{\Big|}
\begin{matrix}
    4& 5\\
    4& 5
    \end{matrix}
\overset{+}{\Big|}
&
\overset{+}{\Big|}
\begin{matrix}
    4& 6\\
    4& 6
    \end{matrix}
\overset{+}{\Big|}
&
\overset{+}{\Big|}
\begin{matrix}
    5& 6\\
    5& 6
    \end{matrix}
\overset{+}{\Big|}
}=
\left[
\begin{matrix}
4&6&12&\\
13&18&27\\
13&18&27\\
40&48&60
\end{matrix}
\right].
\end{equation*}
The matrix $\mathcal{Q}_2(\mathbf C)$ is obtained from  $\mathcal{R}_2(\mathbf C)$ by deleting the  row indexed with $(2,1)$.
\end{example}
\subsection{Links between  matrix ${\mathcal R}_m(\mathbf C)$, matrix $\mathcal B(\mathbf C)$ and symmetrizer}\label{subsection 2.5}

Recall that the matrices $\pi_S(\mathbf T):=(\mathbf T+\mathbf T^T)/2$ and $(\mathbf T-\mathbf T^T)/2$ are called the symmetric part
and skew-symmetric part of a square matrix $\mathbf T$, respectively. The equality
$\mathbf T=(\mathbf T+\mathbf T^T)/2+(\mathbf T-\mathbf T^T)/2$ expresses the well-known fact that
an arbitrary square matrix can be represented
uniquely as a sum of a symmetric matrix and a skew-symmetric matrix.
Similarly, with a general  $m$th-order $K\times \dots\times K$ tensor $\mathcal T$ one can uniquely associate
its symmetric part $\pi_S (\mathcal T)$ --- a tensor whose entry with indices $j_1,\dots,j_m$
is equal to
\begin{equation}
\frac{1}{m!}\sum\limits_{(l_1,\dots,l_m)\in P_{\{j_1,\dots,j_m\}}}  (\mathcal T)_{(l_1,\dots,l_m)}
\label{symmetrization}
\end{equation}
(that is, to get  $\pi_S (\mathcal T)$ we should take the  average  of $m!$ tensors obtained from $\mathcal T$ by
all possible permutations of the indices).
The mapping $\pi_S$ is called symmetrizer (also known as  symmetrization map \cite{marcus1973finite} or
completely symmetric operator \cite{Marcus1964}; in \cite{Schott2003}
a matrix representation of $\pi_S$ was called  Kronecker product permutation matrix).

It is well known that   $m$th-order $K\times\dots \times K$ tensors
can be vectorized into vectors of  $\mathbb R^{K^m}$ in such a way that for any vectors $\mathbf t_1,\dots,\mathbf t_m\in\mathbb R^K$
the rank-1 tensor $\mathbf t_1\circ\dots\circ\mathbf t_m$ corresponds to the vector $\mathbf t_1\otimes\dots\otimes\mathbf t_m$.
This allows us to consider the symmetrizer $\pi_S$ on the space $\mathbb R^{K^m}$. In particular, by
\eqref{symmetrization},
\begin{equation}
\pi_S (\mathbf t_1\otimes\dots\otimes \mathbf t_m)=\frac{1}{m!}\sum\limits_{(l_1,\dots,l_m)\in P_m}\mathbf t_{l_1}\otimes\dots\otimes \mathbf t_{l_m}.\label{eq:deltaKmrank1t}
\end{equation}
The following proposition makes the link between $\mathcal B(\mathbf C)$ and $\mathcal{R}_m(\mathbf C)$ and is the main result of this section.
\begin{proposition}\label{proposition:2.23}
Let $\mathbf C\in\mathbb R^{K\times R}$, $K\leq R$, $m=R-K+2$, and $k_{\mathbf C}\geq K-1$. Let also $\mathcal B(\mathbf C)$ be defined by \eqref{def:B(C)} and let
$\mathcal{R}_m(\mathbf C)^T\upharpoonright_{\textup{range}(\pi_S)}$ denote the
restriction of the mapping $\mathcal{R}_m(\mathbf C)^T:\ \mathbb R^{K^m}\rightarrow \mathbb R^{C^m_R}$ onto $\textup{range}(\pi_S)$. Then
\begin{itemize}
\item[\textup{(i)}] The matrix $\mathcal R_m(\mathbf C)$ has full column rank. Hence, $\dim\textup{range}(\mathcal R_m(\mathbf C)^T)=C^m_R$;
\item[\textup{(ii)}]
$\dim \left( \ker \left(\mathcal{R}_m(\mathbf C)^T\upharpoonright_{\textup{range}(\pi_S)}\right)\right)=C^{K-1}_R$;
\item[\textup{(iii)}] If $k_{\mathbf C} = K$, then
$\ker\left(\mathcal{R}_m(\mathbf C)^T\upharpoonright_{\textup{range}(\pi_S)}\right)= \textup{range} (\mathcal B(\mathbf C)^{(m)})$.
\end{itemize}
\end{proposition}
In the remaining part of this subsection we prove Proposition \ref{proposition:2.23}. 
Readers who are mainly interested in the overall development and algorithms, can safely
skip the rest of this section. We need auxiliary results and notations that we
will also use in Subsection \ref{subsection:3.3}.

Let $\{\mathbf e_j^K\}_{j=1}^K$  denote the canonical basis of $\mathbb R^K$.
Then $\{\mathbf e_{j_1}^K\otimes\dots\otimes \mathbf e_{j_m}^K\}_{(j_1,\dots,j_m)\in R^m_K}$ is the canonical basis of $\mathbb R^{K^m}$
and  by \eqref{eq:deltaKmrank1t}, 
\begin{equation}
\pi_S (\mathbf e_{j_1}^K\otimes\dots\otimes \mathbf e_{j_m}^K)=
\frac{1}{m!}\sum\limits_{(l_1,\dots,l_m)\in P_{\{j_1,\dots,j_m\}}}\mathbf e_{l_1}^K\otimes\dots\otimes \mathbf e_{l_m}^K.\label{eq:deltaKmrank1}
\end{equation}
Let the matrix $\mathbf G\in\mathbb R^{K^m\times C^m_{K+m-1}}$  be defined as follows:
\begin{equation}
\mathbf G \text{ has columns } \{\pi_S(\mathbf e_{j_1}^K\otimes\dots\otimes \mathbf e_{j_m}^K):\ (j_1,\dots,j_m)\in Q^m_K\}.\label{eq:2.11}
\end{equation}
The following lemma follows directly from the definitions of $\pi_S$ and $\mathbf G$ and is well known.
\begin{lemma}\label{lemma:2.14(i)}\cite{Schott2003}
Let $\pi_S$ and $\mathbf G$ be defined by \eqref{eq:deltaKmrank1}--\eqref{eq:2.11}. Then
the columns of the matrix $\mathbf G$ form an orthogonal basis of $\textup{range}(\pi_S)$; in particular,
 $\dim\textup{range}(\pi_S)=C_{K+m-1}^m$.
\end{lemma}

The following lemma explains that the matrix ${\mathcal R}_m(\mathbf C)$ is obtained from $\mathbf C$
by picking all combinations of $m$ columns, and symmetrizing the corresponding rank-1 tensor.
Note that it is the symmetrization that introduces permanents.
\begin{lemma}\label{lemma:permanental compound_equiv}
Let $\mathbf C=\left[\begin{matrix} \mathbf c_1& \dots& \mathbf c_R \end{matrix}\right]\in \mathbb R^{K\times R}$. Then
\begin{equation}\label{RmPCmsym}
\mathcal{R}_m(\mathbf C)=
m!\left[\begin{matrix}
\pi_S(\mathbf c_1\otimes\dots\otimes \mathbf c_m)&
\dots&
\pi_S(\mathbf c_{R-m+1}\otimes\dots\otimes \mathbf c_R)
\end{matrix}\right].
\end{equation}
\end{lemma}
\begin{proof}
By \eqref{eq:deltaKmrank1t},
the $(i_1,\dots,i_m)$-th entry of the vector $m!\pi_S (\mathbf c_{j_1}\otimes\dots\otimes \mathbf c_{j_m})$
is equal to
\begin{equation*}
\begin{split}
\sum\limits_{(l_1,\dots,l_m)\in P_m} c_{i_1j_{l_1}}\cdots c_{i_mj_{l_m}}&=
\textup{perm}
\left[
\begin{matrix}
c_{i_1j_1}&\dots & c_{i_1j_m}\\
\vdots&\vdots&\vdots\\
c_{i_mj_1}&\dots & c_{i_mj_m}
\end{matrix}
\right]
\\
&=\textup{perm\ }\mathbf C((i_1,\dots,i_m),(j_1,\dots,j_m)).
\end{split}
\end{equation*}
Hence,  \eqref{RmPCmsym} follows from Definition \ref{def:extended permanental compound matrices}.
\end{proof}
\begin{example}\label{example:EPC22}
Let the matrix $\mathbf C$ be as in Example \ref{example:EPC2}. Then
$$
\mathcal{R}_2(\mathbf C)^T=2!\left[\begin{matrix}
\frac{1}{2!}\left([1\ 4]\otimes[2\ 5]+[2\ 5]\otimes[1\ 4]\right)\\
\frac{1}{2!}\left([1\ 4]\otimes[3\ 6]+[3\ 6]\otimes[1\ 4]\right)\\
\frac{1}{2!}\left([2\ 5]\otimes[3\ 6]+[3\ 6]\otimes[2\ 5]
\right)
\end{matrix}\right]
=
\left[\begin{matrix}4&13&13&40\\ 6&18&18&48\\ 12&27&27&60\end{matrix}\right].
$$
\end{example}
Let $\left\{\mathbf e_{(j_1,\dots,j_m)}^{C^m_{K+m-1}}\right\}_{(j_1,\dots,j_m)\in Q^m_K}$
denote the canonical basis of $\mathbb R^{C^m_{K+m-1}}$.
Define the $C^m_{K+m-1}$-by-$K^m$  matrix $\mathbf H$ as follows
\begin{equation}
\mathbf H \text{ has columns } \{\mathbf e_{[j_1,\dots,j_m]}^{C^m_{K+m-1}}:\ (j_1,\dots,j_m)\in R^m_K\},\label{eq:2.12}
\end{equation}
in which $[j_1,\dots,j_m]$ denotes the ordered version of
$(j_1,\dots,j_m)$.
For all $K^m$ entries of a symmetric $m$-th order $K\times\dots\times K$ tensor, the corresponding column of $\mathbf H$ contains a ``1''
at the first index combination (in lexicographic ordering) where that entry can be found. The matrix $\mathbf H$ can be used to
``compress'' symmetric $K\times\dots\times K$ tensors by removing redundancies. The matrix $\mathbf G$ above does the opposite thing, so
$\mathbf G$ and $\mathbf H$ act as each other's inverse.
It is easy to prove  that indeed $\mathbf H\mathbf G=\mathbf I_{C^m_{K+m-1}}$.
The relations in the following lemma reflect the same relationship and will be used in Subsection \ref{subsection:3.3}.
\begin{lemma}\label{lemma:2.22}
Let $\mathbf C\in \mathbb R^{K\times R}$ and
let the matrices $\mathbf G$ and $\mathbf H$ be defined
 by \eqref{eq:2.11} and \eqref{eq:2.12}, respectively.
 Then
\begin{itemize}
\item[\textup{(i)}]
$\mathcal R_m(\mathbf C)^T=\mathcal Q_m(\mathbf C)^T\mathbf H$;
\item[\textup{(ii)}]
$\mathcal R_m(\mathbf C)^T\mathbf G= \mathcal Q_m(\mathbf C)^T$.
\end{itemize}
\end{lemma}
\begin{proof}
As the proof is technical, it is given in the supplementary materials.
\end{proof}

{\em Proof of Proposition \ref{proposition:2.23}.}
\textup{(i)}\  Assume that there exists 
$
\widehat{\mathbf t}=[t_{(1,\dots,m)}\ \dots \ t_{(R-m+1,\dots,R)}]^T\in\mathbb R^{C^m_R}
$
such that
$\mathcal R_m(\mathbf C)\widehat{\mathbf t}=\vzero$.
Then, by Lemma \ref{lemma:permanental compound_equiv},
\begin{equation}
\sum\limits_{(p_1,\dots,p_m)\in S^m_R}t_{(p_1,\dots,p_m)}\pi_S(\mathbf c_{p_1}\otimes\dots\otimes \mathbf c_{p_m})=\vzero.\label{eq:2.17}
\end{equation}
Let us fix $(i_1,\dots,i_m)\in S^m_R$ and set
$\{j_1,\dots,j_{K-1}\}:=\{1,\dots,R\}\setminus\{i_1,\dots,i_{m-1}\}$.
Then $i_m\in\{j_1,\dots,j_{K-1}\}$. Without loss of generality we can assume that $j_{K-1}=i_m$.

Since $k_{\mathbf C}\geq K-1$, it follows that there exists a vector $\mathbf y$ such that $\mathbf y$ is orthogonal to the vectors
$\mathbf c_{j_1},\dots,\mathbf c_{j_{K-2}}$, and $\mathbf y$ is not orthogonal to any of $\mathbf c_{i_1},\dots,\mathbf c_{i_{m}}$.
Let $\alpha_{(p_1,\dots,p_m)}$ denote
the $(p_1,\dots,p_m)$-th entry of the vector $\mathcal R_m(\mathbf C)^T(\mathbf y\otimes\dots\otimes\mathbf y)$.
Then, by Lemma \ref{lemma:permanental compound_equiv},
\begin{equation}
\begin{split}
\alpha_{(p_1,\dots,p_m)}=&\pi_S(\mathbf c_{p_1}\otimes\dots\otimes \mathbf c_{p_m})^T(\mathbf y\otimes\dots\otimes\mathbf y)=\\
&\frac{1}{m!}\sum\limits_{(l_1,\dots,l_m)\in P_{\{p_1,\dots,p_m\}}}
(\mathbf c_{l_1}^T\mathbf y)\cdots (\mathbf c_{l_m}^T\mathbf y)=(\mathbf c_{p_1}^T\mathbf y)\cdots (\mathbf c_{p_m}^T\mathbf y).
\end{split}\label{eq:2.9}
\end{equation}
By the construction of $\mathbf y$,  $\alpha_{(p_1,\dots,p_m)}\ne 0$
 if and only if $\{p_1,\dots,p_{m}\}=\{i_1,\dots,i_{m}\}$.
Then, by \eqref{eq:2.17}--\eqref{eq:2.9},
\begin{equation*}
\begin{split}
&0=\sum\limits_{(p_1,\dots,p_m)\in S^m_R}t_{(p_1,\dots,p_m)}\pi_S(\mathbf c_{p_1}\otimes\dots\otimes \mathbf c_{p_m})^T(\mathbf y\otimes\dots\otimes\mathbf y)=\\
&\sum\limits_{(p_1,\dots,p_m)\in S^m_R}t_{(p_1,\dots,p_m)}\alpha_{(p_1,\dots,p_m)}=t_{(i_1,\dots,i_m)}\alpha_{(i_1,\dots,i_m)}.
\end{split}
\end{equation*}
Hence, $t_{(i_1,\dots,i_m)}=0$. Since $(i_1,\dots,i_m)$ was arbitrary we obtain $\widehat{\mathbf t}=\vzero$.

\textup{(ii)}\ From step \textup{(i)}, Lemma \ref{lemma:2.14(i)}, and Lemma \ref{lemma:2.22} \textup{(i)},\textup{(ii)}
it follows that
\begin{equation*}
\begin{split}
C^m_R=&\dim\textup{range}(\mathcal R_m(\mathbf C)^T)\geq
\dim\textup{range}(\mathcal R_m(\mathbf C)^T\upharpoonright_{\textup{range} (\pi_S)})=\\
&\dim\textup{range}(\mathcal R_m(\mathbf C)^T\mathbf G)=
\dim\textup{range}(\mathcal Q_m(\mathbf C)^T)\geq\\
&\dim\textup{range}(\mathcal Q_m(\mathbf C)^T\mathbf H)=
\dim\textup{range}(\mathcal R_m(\mathbf C)^T)=C^m_R.
\end{split}
\end{equation*}
Hence, $\dim\textup{range}(\mathcal R_m(\mathbf C)^T\upharpoonright_{\textup{range}(\pi_S)})=C^m_R$.
By the rank–nullity theorem,
\begin{align*}
\dim\textup{ker}\ (\mathcal R_m(\mathbf C)^T\upharpoonright_{\textup{range} (\pi_S)})=&\dim\textup{range} (\pi_S)-\dim\textup{range} (\mathcal R_m(\mathbf C)^T\upharpoonright_{\textup{range}(\pi_S)})=\\
&C^m_{K+m-1}-C^m_R=C^{R-K+2}_{R+1}-C^{R-K+2}_R=C^{K-1}_R.
\end{align*}

\textup{(iii)}\ 
Let $\mathbf y$ denote the $(j_1,\dots,j_{K-1})$-th column of $\mathcal B(\mathbf C)$. It is clear that the vector $\underbrace{\mathbf y\otimes\dots\otimes \mathbf y}_m$ is contained in $\textup{range}(\pi_S)$.
Hence, $\textup{range}\left(\mathcal B(\mathbf C)^{(m)}\right)\subseteq \textup{range}(\pi_S)$. By step \textup{(ii)} and Proposition \ref{prop:newB(C)} (ii),
$
\dim\textup{ker}\ (\mathcal R_m(\mathbf C)^T\upharpoonright_{\textup{range} (\pi_S)})=C^{K-1}_R=\dim \textup{range}\left( \mathcal B(\mathbf C)^{(m)}\right).
$
To complete the proof we must check that $\mathcal R_m(\mathbf C)^T(\mathbf y\otimes\dots\otimes \mathbf y)=\vzero$ for all $(j_1,\dots,j_{K-1})\in S^{K-1}_R$.
From the construction of the matrix $\mathcal B(\mathbf C)$ it follows that $\mathbf y$ is orthogonal to the vectors $\mathbf c_{j_1},\dots,\mathbf c_{j_{K-1}}$.
Since $(K-1)+m=R+1>R$, it follows that $(\mathbf c_{p_1}^T\mathbf y)\cdots (\mathbf c_{p_m}^T\mathbf y)=0$ for all $(p_1,\dots,p_m)\in S^m_R$.
Hence, by \eqref{eq:2.9}, $\mathcal R_m(\mathbf C)^T(\mathbf y\otimes\dots\otimes \mathbf y)=\vzero$.
\qquad\endproof

The following corollary of Proposition \ref{proposition:2.23} will be used in Subsection \ref{Section:5}.
\begin{corollary}\label{corollary:2.18}
Let the conditions of Proposition \ref{proposition:2.23} hold and  let $k_{\mathbf C}= K-1$. Then
the subspace 
$\ker\left(\mathcal{R}_m(\mathbf C)^T\upharpoonright_{\textup{range}(\pi_S)}\right)$
cannot be spanned by  vectors of the form
$\{\mathbf y_p\otimes \mathbf z_p\}_{p=1}^{C^{K-1}_R}$, where $\mathbf y_p\in\mathbb R^K$ and
$\mathbf z_p\in\mathbb R^{K^{m-1}}$.
\end{corollary}
\begin{proof}
The proof is given in the supplementary materials.
\end{proof}
\section{Transformation of the CPD using polarized compound matrices}\label{sec:multilinear maps}
In this section we derive the crucial expression  \eqref{eq:mainformula}. The matrix  $\RTm{m}$ is constructed from polarized
compound matrices of the slices of the given tensor $\mathcal T$. The entries of polarized compound matrices  are mixed discriminants.
The notions of mixed discriminants and polarized  compound matrices  are introduced in the first two subsections.
\subsection{Mixed discriminants}
The mixed discriminant is variant of the determinant that has more than one matrix argument.
\begin{definition}\cite{Aleksandrov1938}\label{def:mixeddet}
Let  $\mathbf T_1,\dots,\mathbf T_m\in\mathbb R^{m\times m}$.
The mixed discriminant,  denoted by $\mathcal D(\mathbf T_1,\dots,\mathbf T_m)$, is defined as
the coefficient of $x_1\cdots x_m$ in $\det (x_1\mathbf T_1+\dots+x_m\mathbf T_m)$, that is
\begin{equation}\label{eq:firstrepresD}
\mathcal D(\mathbf T_1,\dots,\mathbf T_m)=
\left.
\frac{\partial^m\left( \det (x_1\mathbf T_1+\dots+x_m\mathbf T_m)\right)}{\partial x_1\dots\partial x_m}\right|_{x_1=\dots=x_m=0}.
\end{equation}
\end{definition}
For convenience, we have dropped the factor $1/m!$ before the fraction in
\eqref{eq:firstrepresD}.
Definition \ref{def:mixeddet}   implies the following lemmas.
\begin{lemma}\cite{Aleksandrov1938}\label{lemma:multidet}
The mapping $(\mathbf T_1,\dots,\mathbf T_m)\rightarrow \mathcal D(\mathbf T_1,\dots,\mathbf T_m)$ is multilinear and symmetric in its arguments.
\end{lemma}
\begin{lemma} \cite{Egorychev1981}\label{lemma:Egorychev} Let $\mathbf d_1,\dots,\mathbf d_m\in\mathbb R^m$. Then
$
\mathcal D\left(\textup{\text{Diag}}(\mathbf d_1),\dots, \textup{\text{Diag}}(\mathbf d_m)\right)=
\textup{perm\ }\left[\begin{matrix}\mathbf d_1&\dots&\mathbf d_m\end{matrix}\right].
$
\end{lemma}

{\em Proof.}
\begin{equation*}
\begin{split}
&\mathcal D (\textup{\text{Diag}}(\left[\begin{matrix}d_{11}&\dots &d_{m1}\end{matrix}\right]),\dots,\textup{\text{Diag}}(\left[\begin{matrix}d_{1m}&\dots & d_{mm}\end{matrix}\right]))=\\
&
\left.
\frac{\partial^m\left((x_1d_{11}+\dots+ x_md_{1m})\cdots (x_1d_{m1}+\dots+ x_md_{mm})\right)}{\partial x_1\dots\partial x_m}\right|_{x_1=\dots=x_m=0}=\\
&\sum\limits_{(l_1,\dots,l_m)\in P_m} d_{1l_1}\cdots d_{ml_m}=
\textup{perm\ }\left[\begin{matrix}\mathbf d_1&\dots&\mathbf d_m\end{matrix}\right].\qquad\endproof
\end{split}
\end{equation*}
Mixed discriminants may be computed numerically from \eqref{eq:firstrepresD}. A direct expression in terms of determinants is given in the following lemma.
\begin{lemma}\cite{determinantfunctions, Bapat1989}\label{lemma:equivdet}
Let  $\mathbf T_1,\dots,\mathbf T_m\in\mathbb R^{m\times m}$. Then
\begin{equation}\label{eq:repwithpolar}
\mathcal D(\mathbf T_1,\dots,\mathbf T_m)=
\sum\limits_{k=1}^m(-1)^{m-k}
\sum\limits_{1\leq i_1<i_2<\dots<i_k\leq m}\det(\mathbf T_{i_1}+\dots+\mathbf T_{i_k}).
\end{equation}
\end{lemma}
The way in which \eqref{eq:repwithpolar} obtains the mixed discriminant from the determinant is an instance of a technique  called polarization \cite{Landsberg}.
\subsection{Polarized compound matrices}
Let $m\geq 2$. In this subsection we discuss a polarized version of compound matrices, in which the mixed discriminant replaces the determinant.
\begin{definition}\label{def:3.4}
Let $\min(I,J)\geq m\geq 2$ and let  $\mathbf T_1,\dots,\mathbf T_m\in\mathbb R^{I\times J}$.
The $C_I^m$-by-$C_J^m$ matrix
$\rdm{m-1}(\mathbf T_1,\dots,\mathbf T_m)$ is defined by
\begin{equation}\label{eq:firstrepresDcomp}
\rdm{m-1}(\mathbf T_1,\dots,\mathbf T_m)=
\left.
\frac{\partial^m\left( \mathcal C_m (x_1\mathbf T_1+\dots+x_m\mathbf T_m)\right)}{\partial x_1\dots\partial x_m}\right|_{x_1=\dots=x_m=0}.
\end{equation}
\end{definition}

In the following lemmas we establish properties of $\rdm{m-1}(\mathbf T_1,\dots,\mathbf T_m)$.

\begin{lemma}\label{lemma:rankdetext}
Let $\mathbf T\in\mathbb R^{I\times J}$ and $\mathbf d_1,\dots,\mathbf d_m\in\mathbb R^R$. Then
\begin{itemize}
\item[\textup{(i)}] the mapping $(\mathbf T_1,\dots,\mathbf T_m)\rightarrow \rdm{m-1}(\mathbf T_1,\dots,\mathbf T_m)$ is multilinear and symmetric in its arguments;
\item[\textup{(ii)}]
an equivalent expression
for $\rdm{m-1}(\mathbf T_1,\dots,\mathbf T_m)$ is
\begin{equation*}
\rdm{m-1}(\mathbf T_1,\dots,\mathbf T_m)=
\sum\limits_{k=1}^m(-1)^{m-k}
\sum\limits_{1\leq i_1<i_2<\dots<i_k\leq m}\mathcal C_m(\mathbf T_{i_1}+\dots+\mathbf T_{i_k});
\end{equation*}
\item[\textup{(iii)}]  $\rdm{m-1}(\mathbf T,\dots, \mathbf T)=m!\mathcal C_m(\mathbf T)$;
\item[\textup{(iv)}]  $r_{\mathbf T}\leq m-1$ if and only if $\rdm{m-1}(\mathbf T,\dots, \mathbf T)=\mzero$;
\item[\textup{(v)}]
 $
 \rdm{m-1}\left(\textup{\text{Diag}}(\mathbf d_1),\dots, \textup{\text{Diag}}(\mathbf d_m)\right)=
\textup{\text{Diag}}\left(\mathcal{PC}_m(\left[\begin{matrix}\mathbf d_1&\dots&\mathbf d_m\end{matrix}\right])\right).
$
\end{itemize}
\end{lemma}
\begin{proof}
From Definitions \ref{def:compound matrices} and \ref{def:3.4} it follows that the $(i, j)$-th entry of the matrix   $\rdm{m-1}(\mathbf T_1,\dots,\mathbf T_m)$ is equal to
$\mathcal D(\mathbf T_1(S_I^m(i), S_J^m(j)),\dots,\mathbf T_m(S_I^m(i), S_J^m(j)))$.
Hence,  statements \textup{(i)} and \textup{(ii)} follow from
 Lemma \ref{lemma:multidet} and Lemma \ref{lemma:equivdet}, respectively.
Statement \textup{(iii)} follows from \eqref{eq:firstrepresDcomp}. Statement \textup{(iv)} follows from \textup{(iii)} and
Lemma \ref{compoundprop1} (2). Finally, \textup{(v)} follows from Lemma \ref{Lemma2.3}, statement \textup{(ii)}, and
Lemma \ref{lemma:Egorychev}.
\end{proof}
\begin{example}
\begin{equation}
\begin{split}
&\rdm{2}(\mathbf T_1,\mathbf T_2,\mathbf T_3)= \mathcal C_3(\mathbf T_1+\mathbf T_2+\mathbf T_3)-\\
&\qquad \mathcal C_3(\mathbf T_1+\mathbf T_2)-\mathcal C_3(\mathbf T_1+\mathbf T_3)-\mathcal C_3(\mathbf T_2+\mathbf T_3)+
\mathcal C_3(\mathbf T_1)+\mathcal C_3(\mathbf T_2)+\mathcal C_3(\mathbf T_3).
\end{split}\label{eq:example123}
\end{equation}
\end{example}
\begin{remark}
 The polarized compound matrix is a matrix representation of the higher-order tensor obtained by the low-rank detection mapping in
\cite{DeLathauwer2006, DimLievenLL1}. More specifically, in \cite{DeLathauwer2006} a rank-$1$ detection mapping ($m=2$) was used to compute the CPD and in
\cite{DimLievenLL1} a rank-($L,L,1$) detection mapping ($m$ arbitrary) was used to compute the decomposition in rank-($L,L,1$) terms. Statement \textup{(iv)} of Lemma \ref{lemma:rankdetext} explains the terminology.
\end{remark}

The following
counterpart of Lemma \ref{LemmaCompound} holds for polarized compound matrices.
\begin{lemma}\label{lemma:3.7}
Let $\mathbf A\in \mathbb R^{I\times R}$, $\mathbf B\in \mathbb R^{J\times R}$, $\mathbf d_1,\dots,\mathbf d_m\in\mathbb R^R$, and $m\leq\min(I,J,R)$. Then
\begin{equation}\label{eq:3.4}
\begin{split}
\rdm{m-1}&\left(\mathbf A\textup{\text{Diag}}(\mathbf d_1)\mathbf B^T,\dots, \mathbf A\textup{\text{Diag}}(\mathbf d_m)\mathbf B^T\right)=\\
&\mathcal C_m(\mathbf A)\textup{\text{Diag}}\left(\mathcal{PC}_m(\left[\begin{matrix}\mathbf d_1&\dots&\mathbf d_m\end{matrix}\right])\right)\mathcal C_m(\mathbf B)^T.
\end{split}
\end{equation}
\end{lemma}
\begin{proof}
From Lemma \ref{lemma:rankdetext} \textup{(ii)} and Lemma \ref{LemmaCompound}  we have
\begin{equation}\label{eq:3.7}
\begin{split}
\rdm{m-1}&\left(\mathbf A\textup{\text{Diag}}(\mathbf d_1)\mathbf B^T,\dots, \mathbf A\textup{\text{Diag}}(\mathbf d_m)\mathbf B^T\right)=\\
&\mathcal C_m(\mathbf A)
\rdm{m-1}\left(\textup{\text{Diag}}(\mathbf d_1),\dots, \textup{\text{Diag}}(\mathbf d_m)\right)
\mathcal C_m(\mathbf B)^T.
\end{split}
\end{equation}
Now \eqref{eq:3.4} follows from \eqref{eq:3.7} and Lemma \ref{lemma:rankdetext} \textup{(v)}.
\end{proof}
\subsection{Transformation of the tensor}\label{subsection:3.3}
We stack polarized compound matrices obtained from the slices  of a given tensor in matrices $\RTm{m}$ and $\TTm{m}$.
In $\RTm{m}$ we consider all slice combinations, while in $\TTm{m}$ we avoid doubles by taking into account the invariance of polarized compound matrices
under permutation of their arguments. In our algorithms we will work with the smaller matrix $\TTm{m}$ while in the theoretical
development we will use $\RTm{m}$.
\begin{definition}\label{definition:detecting matrix}
Let $\mathcal T$ be an $I\times J\times K$ tensor with frontal slices $\mathbf T_1,\dots,\mathbf T_K\in\mathbb R^{I\times J}$. The
$(j_1,\dots,j_m)$-th column of the
$C^m_IC^m_J$-by-$K^m$ (resp. $C^m_IC^m_J$-by-$C^m_{K+m-1}$) matrix $\RTm{m}$ (resp. $\TTm{m}$) equals $\textup{vec}\left( \rdm{m-1}(\mathbf T_{j_1},\dots,\mathbf T_{j_m})\right)$,
where $(j_1,\dots,j_m)\in R^m_K$ (resp. $Q^m_K$).
\end{definition}

Let $\RTm{m}\upharpoonright_{\textup{range}(\pi_S)}$ denote the
restriction of the mapping $\RTm{m}:\ \mathbb R^{K^m}\rightarrow \mathbb R^{C^m_IC^m_J}$ onto $\textup{range}(\pi_S)$.
In the following lemma we express the matrices $\RTm{m}$ and $\TTm{m}$ via the factor matrices of $\mathcal T$ and make a link between the kernel
of $\RTm{m}\upharpoonright_{\textup{range}(\pi_S)}$ and $\TTm{m}$. These results are key to our overall derivation.
\begin{lemma}\label{lemma:3.9}
Let  $\mathbf A\in \mathbb R^{I\times R}$, $\mathbf B\in \mathbb R^{J\times R}$,  $\mathbf C\in \mathbb R^{K\times R}$, and $\mathcal T=[\mathbf A,\mathbf B,\mathbf C]_R$. Then for $m\leq\min(I,J,K,R)$,
\begin{itemize}
\item[\textup{(i)}]
$\RTm{m}=\left[\mathcal C_m(\mathbf A)\odot\mathcal C_m(\mathbf B)\right]\mathcal R_m(\mathbf C)^T$;
\item[\textup{(ii)}]
$\TTm{m}=\left[\mathcal C_m(\mathbf A)\odot\mathcal C_m(\mathbf B)\right]\mathcal Q_m(\mathbf C)^T$;
\item[\textup{(iii)}]
$\ker( \RTm{m}\upharpoonright_{\textup{range}(\pi_S)})=\mathbf G\ker(
\TTm{m})$, where  $\mathbf G$ is defined in \eqref{eq:2.11}.
\end{itemize}
\end{lemma}
\begin{proof}
\textup{(i)}\
Let $\mathbf c^{1},\dots,\mathbf c^{K}$ be the columns of the matrix $\mathbf C^T$.
Recall that the frontal slices of $\mathcal T$ can be expressed as in  \eqref{eq:matranalog1}.
Then, by
Lemma \ref{lemma:3.7} and identity \eqref{eqmatrtovec},
\begin{equation*}
\begin{split}
\textup{vec}&\left(\rdm{m-1} \left(\mathbf A\textup{\text{Diag}}(\mathbf c^{j_1})\mathbf B^T,\dots, \mathbf A\textup{\text{Diag}}(\mathbf c^{j_m})\mathbf B^T\right)       \right)=\\
\textup{vec}&\left(
\mathcal C_m(\mathbf A)\textup{\text{Diag}}\left(\mathcal{PC}_m(\left[\begin{matrix}\mathbf c^{j_1}&\dots&\mathbf c^{j_m}\end{matrix}\right])\right)\mathcal C_m(\mathbf B)^T\right)=\\
&\left[\mathcal C_m(\mathbf A)\odot\mathcal C_m(\mathbf B)\right]\mathcal{PC}_m(\left[\begin{matrix}\mathbf c^{j_1}&\dots&\mathbf c^{j_m}\end{matrix}\right]).
\end{split}
\end{equation*}
Now \textup{(i)} and  \textup{(ii)} follow from Definition \ref{definition:detecting matrix} and Lemma \ref{lemma:permanental compound_equiv2}.

\textup{(iii)}\
 From \textup{(i)}, \textup{(ii)}, and Lemma \ref{lemma:2.22} \textup{(ii)}
 it follows that $\RTm{m}\mathbf G=\TTm{m}$. Since, by Lemma \ref{lemma:2.14(i)}, $\textup{range}(\pi_S)=\textup{range}(\mathbf G)$ we obtain \textup{(iii)}.
 \end{proof}
\section{Overall results and algorithms}\label{Section4}
\subsection{Algorithm \ref{alg:Main}}
Overall, Algorithm \ref{alg:Main} goes now as follows. We first compute
$\TTm{m}$ from $\mathcal T$, determine its null space, which, after symmetrization, yields
$\ker\left(\RTm{m}\upharpoonright_{\textup{range}(\pi_S)}\right)$,
as explained in Lemma \ref{lemma:3.9} \textup{(iii)}. The following lemma makes now, for a particular choice of $m$, a connection with $\mathcal B(\mathbf C)$.
\begin{lemma}\label{lemma 4.1iandii}
Let   $\mathcal T=[\mathbf A,\mathbf B,\mathbf C]_R$,  $m:= R-K+2$, and
let $ \RTm{m}$ be defined in Definition \ref{definition:detecting matrix}.
Assume that
$k_{\mathbf C}=K$ and  that
$\mathcal C_m(\mathbf A)\odot \mathcal C_m(\mathbf B)$  has full column rank.
Then
\begin{itemize}
\item[\textup{(i)}] $\ker(\RTm{m}\upharpoonright_{\textup{range}(\pi_S)})=\textup{range}(\mathcal B(\mathbf C)^{(m)})$;
\item[\textup{(ii)}]     $\dim\ker(\RTm{m}\upharpoonright_{\textup{range}(\pi_S)})=C^{K-1}_R$;
\end{itemize}
\end{lemma}
\begin{proof}
Since $\mathcal C_m(\mathbf A)\odot \mathcal C_m(\mathbf B)$  has full column rank, it follows from Lemma \ref{lemma:3.9} \textup{(i)} that
$\ker(\RTm{m}\upharpoonright_{\textup{range}(\pi_S)})=
\ker\left(\mathcal{R}_m(\mathbf C)^T\upharpoonright_{\textup{range}(\pi_S)}\right)$.
Statements \textup{(i)} and \textup{(ii)} now follow from Proposition \ref{proposition:2.23} \textup{(iii)} and \textup{(ii)}, respectively.
\end{proof}

So far, we have obtained from $\mathcal T$  a basis for the column space of $\mathcal B(\mathbf C)^{(m)}$.
The following lemma explains that the basis vectors may be stacked in a tensor that has
$\mathcal B(\mathbf C)$  as factor matrix. Moreover, the CPD may be computed by a GEVD as in Theorem
 \ref{theoremoralggeneigdec}.
\begin{lemma}\label{lemma 4.1}
Suppose that  the conditions of Lemma \ref{lemma 4.1iandii} hold.
Let $\mathbf W$ be a $K^m\times C^{K-1}_R$ matrix such that
\begin{equation}\label{eq:wker}
\ker(\RTm{m}\upharpoonright_{\textup{range}(\pi_S)})=\textup{range }(\mathbf W)
\end{equation}
 and let
$\mathcal W$ be the $K\times K^{m-1}\times C^{K-1}_R$ tensor such that $\mathbf W=\Matr{\mathcal W}$.
Then
\begin{itemize}
\item[\textup{(i)}]  there exists a nonsingular $C^{K-1}_R\times C^{K-1}_R$
matrix $\mathbf M$ such that
\begin{equation}\label{eq:newCPD}
\mathcal W=\big[\mathcal B(\mathbf C),\mathcal B(\mathbf C)^{(m-1)},\mathbf M\big]_{C^{K-1}_R};
\end{equation}
\item[\textup{(ii)}] $r_{\mathcal W}=C^{K-1}_R$ and the CPD of $\mathcal W$ is unique and can be found algebraically.
\end{itemize}
\end{lemma}
\begin{proof}
\textup{(i)}\quad From Lemma \ref{lemma 4.1iandii} \textup{(ii)} and \eqref{eq:wker} it follows that there exists
a nonsingular $C^{K-1}_R\times C^{K-1}_R$ matrix $\mathbf M$, such that
$
\mathbf W=\mathcal B(\mathbf C)^{(m)}\mathbf M^T=
\left(\mathcal B(\mathbf C)\odot \mathcal B(\mathbf C)^{(m-1)}\right)\mathbf M^T.
$
Hence, by \eqref{eqT_V},  \eqref{eq:newCPD} holds.

\textup{(ii)}\quad From Proposition \ref{prop:newB(C)} it follows that $k_{\mathcal B(\mathbf C)}\geq 2$ and that
 the matrix $\mathcal B(\mathbf C)^{(m-1)}$ has rank $C^{K-1}_R$.
The statement now follows from Theorem \ref{theoremoralggeneigdec}.
\end{proof}

After finding $\mathcal B(\mathbf C)$ up to column permutation and scaling,
we may find $\mathbf C$ as explained in Subsection \ref{subsection 2.1}.
The following Lemma completes the proof of Theorem \ref{theorem: main} \textup{(ii)}.
Its proof shows how the other factor matrices may be determined once $\mathbf C$
 has been obtained. The computation involves another CPD of the form in Theorem \ref{theoremoralggeneigdec}.
The result is a variant of
\cite[Theorem 3.8]{MikaelHO}; in this step of the derivation we do not assume that the decomposition is canonical.

\begin{lemma}\label{lemma:4.2}
Let  $\mathcal T=[\mathbf A,\mathbf B,\mathbf C]_R$ and the $K\times R$ matrix $\mathbf C$ be known.
Assume that
$k_{\mathbf C}=K\geq 2$, and  that $\min(k_{\mathbf A},k_{\mathbf B})+k_{\mathbf C}\geq R+2$.
Then the matrices $\mathbf A$, $\mathbf B$ can be found algebraically up to column scaling.
\end{lemma}
\begin{proof}
We obviously have  $k_{\mathbf C}=r_{\mathbf C}=K$.
Let $\mathbf X=\left[\begin{matrix}\mathbf c_1&\dots&\mathbf c_K\end{matrix}\right]$.
By multiplying with $\mathbf X^{-1}$ we will create a CPD with $K-2$ terms less than $R$.
It is clear that the matrix formed by the first two rows of ${\mathbf X}^{-1}\mathbf C$  has the form $\left[\begin{matrix}\mathbf I_2&\mzero_{2\times (K-2)}&\mathbf Y\end{matrix}\right]$, where $\mathbf Y$ is a $2\times(R-K)$ matrix. Define
$
\widetilde{\mathbf A}:=\left[\begin{matrix}\mathbf a_1&\mathbf a_2&\mathbf a_{K+1}&\dots&\mathbf a_R\end{matrix}\right]$,
$
\widetilde{\mathbf B}:=\left[\begin{matrix}\mathbf b_1&\mathbf b_2&\mathbf b_{K+1}&\dots&\mathbf b_R\end{matrix}\right]$, and
$
\widetilde{\mathbf C}=\left[\begin{matrix}\mathbf I_2&\mathbf Y\end{matrix}\right].
$
Let also $\widetilde{\mathcal T}$ denote the $I\times J\times 2$ tensor such that  $\Matr{\widetilde{\mathcal T}}$ coincides with the first two columns of the matrix
$\Matr{\mathcal T}{\mathbf X}^{-T}$.
From
\eqref{eqT_V} it follows that $\Matr{\mathcal T}{\mathbf X}^{-T}=(\mathbf A\odot \mathbf B)\mathbf C^T{\mathbf X}^{-T}$.
Hence, $\Matr{\widetilde{\mathcal T}}=(\widetilde{\mathbf A}\odot \widetilde{\mathbf B})\widetilde{\mathbf C}^T$ or
$\widetilde{\mathcal T}=[\widetilde{\mathbf A},\widetilde{\mathbf B},\widetilde{\mathbf C}]_{R-K+2}$, which is of the desired form.

It is easy to show that $\widetilde{\mathcal T}$ satisfies the conditions of Theorem \ref{theoremoralggeneigdec}, which means that its rank is $R-K+2$,
that its CPD is unique and that the factor matrices may be found algebraically. The indeterminacies in
$\widetilde{\mathcal T}=[\widehat{\mathbf A},\widehat{\mathbf B},\widehat{\mathbf C}]_{R-K+2}$ are limited to the existence of
a   permutation matrix $\mathbf P$   and  a nonsingular
diagonal matrix $\mathbf{\Lambda}$ such that
$\widetilde{\mathbf C}=\widehat{\mathbf C}\mathbf P\mathbf{\Lambda}$ and
$\widetilde{\mathbf A}\odot \widetilde{\mathbf B}=(\widehat{\mathbf A}\odot \widehat{\mathbf B})\mathbf P\mathbf{\Lambda}^{-1}$.

So far we have  algebraically found the columns of the matrices $\mathbf A$, $\mathbf B$, and hence $\mathbf A\odot\mathbf B$,
with indices in $I:=\{1,2,K+1,\dots,R\}$.
Let $\bar{\mathbf A}$, $\bar{\mathbf B}$, and $\bar{\mathbf C}$ be the submatrices of $\mathbf A$, $\mathbf B$ and $\mathbf C$, respectively,
formed by the columns with indices in $\{3,\dots,K\}$. We now subtract the rank-1 terms that we already know to obtain
$
\mathcal T-\sum\limits_{r\in I}\mathbf a_r\circ\mathbf b_r\circ\mathbf c_r=
\left[\bar{\mathbf A}, \bar{\mathbf B},\bar{\mathbf C} \right]_{K-2}=:\bar{\mathcal T}\  \text{ or }\
(\bar{\mathbf A}\odot \bar{\mathbf B})\bar{\mathbf C}^T=\Matr{\bar{\mathcal T}}.
$
Since the matrix $\bar{\mathbf C}$ has full column rank, the columns of the matrix $\mathbf A\odot\mathbf B$ with indices in  $\{3,\dots,K\}$ coincide with the columns of $\Matr{\bar{\mathcal T}}\bar{\mathbf C}^\dagger$.
Now that also the columns of $\mathbf A\odot\mathbf B$  with indices in
$\{3,\dots,K\}$ have been found, $\mathbf a_r$ and $\mathbf b_r$ are easily obtained by understanding that $\mathbf a_r\otimes\mathbf b_r=\vect{\mathbf b_r\mathbf a_r^T}$, $r=1,\dots,R$.
\end{proof}

\begin{algorithm}
\caption{(Computation of $\mathbf C$, then $\mathbf A$ and $\mathbf B$)}
\label{alg:Main}
\begin{algorithmic}[1]
\Require $\mathcal T\in\mathbb R^{I\times J\times K}$ and $R\geq 2$  with the property that  there exist $\mathbf A\in\mathbb R^{I\times R}$, $\mathbf B\in\mathbb R^{J\times R}$, and $\mathbf C\in\mathbb R^{K\times R}$ such that $\mathcal T=[\mathbf A,\mathbf B,\mathbf C]_R$, $k_{\mathbf C}=K\geq 2$, and
$\mathcal C_m(\mathbf A)\odot \mathcal C_m(\mathbf B)$  has full column rank for $m=R-K+2$.
\Ensure Matrices $\mathbf A\in\mathbb R^{I\times R}$, $\mathbf B\in\mathbb R^{J\times R}$ and $\mathbf C\in\mathbb R^{K\times R}$ such that $\mathcal T=[\mathbf A,\mathbf B,\mathbf C]_R$
\Algphase{{\bf Phase 1 (based on Lemma \ref{lemma 4.1}):} Find the matrix $\mathbf F\in\mathbb R^{K\times C^{K-1}_R}$ such that $\mathbf F$ coincides with $\mathcal B(\mathbf C)$ up to (unknown) column permutation and  scaling}
\Statex {\em Apply Lemma \ref{lemma:3.9} (iii) to find the} $K^m\times C^{K-1}_R$ {\em matrix } $\mathbf W$ {\em such that \eqref{eq:wker} holds}
\State Construct  the $C^m_IC^m_J$-by-$C^m_{K+m-1}$ matrix $\TTm{m}$ by Definition \ref{definition:detecting matrix}
\State Find $\bar{\mathbf w}_1,\dots,\bar{\mathbf w}_{C^{K-1}_R}$ which form a basis  of \textup{ker}($\TTm{m}$)
\State $\mathbf W\leftarrow \mathbf G\left[\begin{matrix}\bar{\mathbf w}_1&\dots&\bar{\mathbf w}_{C^{K-1}_R}\end{matrix}\right]$, where
 $\mathbf G\in\mathbb R^{K^m\times C^m_{K+m-1}}$  is defined by \eqref{eq:2.11}
\Statex {\em Apply Theorem \ref{theoremoralggeneigdec} \textup{(ii)} to find $\mathbf F$}
\State ${\mathcal W} \leftarrow \UnMatr{\mathbf W,K,K^{m-1}}$
\State Compute the CPD $\mathcal W=[\mathbf F,\mathbf F_2,\mathbf F_3]_{C^{K-1}_R}$ ({\em $\mathbf F_2$ and $\mathbf F_3$ are a by-product})\ \ (GEVD)
\Algphase{{\bf Phase 2 (based on properties \eqref{P3}--\eqref{P4}):} Find the matrix $\mathbf C$}
\State Compute $R$ subsets of $C^{K-2}_{R-1}$ columns of $\mathbf F$ that are linearly dependent
\State Compute $\mathbf c_1,\dots,\mathbf c_R$ as orthogonal complements to sets found in step 6
\Algphase{{\bf Phase 3 (based on Lemma \ref{lemma:4.2}):} Find the matrices $\mathbf A$ and $\mathbf B$}
\Statex{\em Apply Theorem \ref{theoremoralggeneigdec} \textup{(ii)} to find the columns of $\mathbf S:=\mathbf A\odot\mathbf B$ with indices in $\{1,2,K+1,\dots,R\}$}
\State $\mathbf Z=\left[\begin{matrix}\mathbf z_1&\dots&\mathbf z_K\end{matrix}\right] \leftarrow \Matr{\mathcal T}\left[\begin{matrix}\mathbf c_1&\dots&\mathbf c_K\end{matrix}\right]^{-T}$
\State $\widetilde{\mathcal T} \leftarrow \UnMatr{\left[\begin{matrix}\mathbf z_1&\mathbf z_2\end{matrix}\right],I,J}$
\State Compute the CPD $\widetilde{\mathcal T}=[\widehat{\mathbf A},\widehat{\mathbf B},\widehat{\mathbf C}]_{R-K+2}$
\qquad\qquad\qquad\qquad\qquad\qquad\quad(GEVD)
\State $\widetilde{\mathbf C}\leftarrow\left[\begin{matrix}\mathbf I_2&\mathbf Y\end{matrix}\right]$, where
$\mathbf Y$ is the $2\times R$ submatrix in upper right-hand corner  of $\left[\begin{matrix}\mathbf c_1&\dots&\mathbf c_K\end{matrix}\right]^{-1}\mathbf C$
\State Compute  permutation matrix $\mathbf P$   and diagonal matrix $\mathbf{\Lambda}$ such that
$\widetilde{\mathbf C}=\widehat{\mathbf C}\mathbf P\mathbf{\Lambda}$
\State $\left[\begin{matrix}\mathbf s_1&\mathbf s_2&\mathbf s_{K+1}\dots&\mathbf s_R\end{matrix}\right]\leftarrow (\widehat{\mathbf A}\odot \widehat{\mathbf B})\mathbf P\mathbf{\Lambda}^{-1}$
\Statex {\em Find the columns of $\mathbf S$ with indices in $\{3,\dots,K\}$}
\State $\mathbf Z \leftarrow \Matr{\mathcal T}-(\widetilde{\mathbf A}\odot\widetilde{\mathbf B})
\left[\begin{matrix}\mathbf c_1&\mathbf c_2&\mathbf c_{K+1}&\dots&\mathbf c_R\end{matrix}\right]^T$
\State  $\left[\begin{matrix}\mathbf s_3&\dots&\mathbf s_K\end{matrix}\right]\leftarrow \mathbf Z\left[\begin{matrix}\mathbf c_3&\dots&\mathbf c_K\end{matrix}\right]^\dagger$
\State Find the columns of $\mathbf A$ and $\mathbf B$ from the equations
$\mathbf a_r\otimes\mathbf b_r=\mathbf s_r$, $r=1,\dots,R$
\end{algorithmic}
\end{algorithm}
\begin{algorithm}
\caption{(Computation of $\mathbf A$ and $\mathbf B$, then $\mathbf C$)}
\label{alg:Main2}
\begin{algorithmic}[1]
 \Require $\mathcal T\in\mathbb R^{I\times J\times K}$ and $R\geq 2$  with the property that  there exist $\mathbf A\in\mathbb R^{I\times R}$, $\mathbf B\in\mathbb R^{J\times R}$, and $\mathbf C\in\mathbb R^{K\times R}$ such that $\mathcal T=[\mathbf A,\mathbf B,\mathbf C]_R$, $k_{\mathbf C}=K\geq 2$, and
$\mathcal C_m(\mathbf A)\odot \mathcal C_m(\mathbf B)$  has full column rank for $m=R-K+2$.
\Ensure Matrices $\mathbf A\in\mathbb R^{I\times R}$, $\mathbf B\in\mathbb R^{J\times R}$ and $\mathbf C\in\mathbb R^{K\times R}$ such that $\mathcal T=[\mathbf A,\mathbf B,\mathbf C]_R$
\Algphase{{\bf Phase 1 (based on Lemma \ref{lemma 4.1}):} Find the matrix $\mathbf F\in\mathbb R^{K\times C^{K-1}_R}$ such that $\mathbf F$ coincides with
$\mathcal B(\mathbf C)$ up to (unknown) column permutation and  scaling}
\Statex \hspace{-5.5ex} {\footnotesize 1--5:}  Identical to Algorithm \ref{alg:Main}
\Algphase{{\bf Phase 2 (based on Lemma \ref{lemma:4.3}):} Find the matrices $\mathbf A$ and $\mathbf B$}
\setalglineno{6}
\State $\mathcal V\leftarrow\UnMatr{\mathbf V,I,J}$, where $\mathbf V=\Matr{\mathcal T}\mathbf F$
\Statex {\em Let $\mathbf V_1,\dots,\mathbf V_{C^{K-1}_R}$ denote the frontal slices of $\mathcal V$}
\Statex {\em Let $\mathcal V_{ij}$ denote the  tensor with frontal slices $\mathbf V_i$ and $\mathbf V_j$}
\Statex {\em Find $C^m_RC^2_m$ pairs $(i,j)$ such that $r_{{\mathcal V}_{ij}}=m$}
\State
$
\mathcal J\leftarrow\{(i,j):\ r_{[\mathbf V_i\ \mathbf V_j]}=r_{[\mathbf V_i^T\ \mathbf V_j^T]}=m,\ 1\leq i<j\leq C^{K-1}_R\}
$
\Statex{\em Apply Theorem \ref{theoremoralggeneigdec} \textup{(ii)} to find  $C^m_RC^2_m$ sets of $m$ columns of $\mathbf A$ and $\mathbf B$ each}
\State Find the CPD  $\mathcal V_{ij}=[{\mathbf A}_{ij},{\mathbf B}_{ij},{\mathbf C}_{ij}]_m$ for each $(i,j)\in \mathcal J$
\Statex ({\em ${\mathbf C}_{ij}$ are a by-product})
\qquad\qquad\qquad\qquad  \hfill (GEVD)
\State $\widetilde{\mathbf A}\leftarrow $ the $(I\times mC^m_RC^2_m)$ matrix formed by the columns of the matrices ${\mathbf A}_{ij}$
\State $\widetilde{\mathbf B}\leftarrow $ the $(I\times mC^m_RC^2_m)$ matrix formed by the columns of the matrices ${\mathbf B}_{ij}$
\State Choose $r_1,\dots,r_R$ such that the sets
$\{\widetilde{\mathbf a}_{r_1},\dots,\widetilde{\mathbf a}_{r_R}\}$ and
$\{\widetilde{\mathbf b}_{r_1},\dots,\widetilde{\mathbf b}_{r_R}\}$ do not
contain collinear vectors
\State $\mathbf A\leftarrow [\begin{matrix}\widetilde{\mathbf a}_{r_1}&\dots&\widetilde{\mathbf a}_{r_R}\end{matrix}]$,
 $\mathbf B\leftarrow [\begin{matrix}\widetilde{\mathbf b}_{r_1}&\dots&\widetilde{\mathbf b}_{r_R}\end{matrix}]$
\Algphase{{\bf Phase 3:} Find the matrix $\mathbf C$}
\State $\mathbf C\leftarrow\left((\mathbf A\odot\mathbf B)^\dagger\Matr{\mathcal T}\right)^T$
\end{algorithmic}
\end{algorithm}
The overall procedure that constitutes the proof of Theorem \ref{theorem: main} (ii) is summarized in Algorithm \ref{alg:Main}.
Phase 2 is formulated in a way that has combinatorial complexity and quickly becomes computationally infeasible. The amount of work may be reduced by exploiting the dependencies in $\mathbf F$ only partially.
\subsection{Algorithm \ref{alg:Main2}}
We derive an algorithmic  variant that further reduces the computational cost.
This algorithm is given in Algorithm \ref{alg:Main2} below.
While Algorithm \ref{alg:Main} first determines $\mathbf C$ and then finds $\mathbf A$ and $\mathbf B$, Algorithm   \ref{alg:Main2} works the other way around.
The basic idea is as follows. Like in Algorithm \ref{alg:Main}, we first find a matrix $\mathbf F$ that is equal to $\mathcal B(\mathbf C)$
up to column permutation and scaling. If $\mathbf C$ is square, we have from Subsection \ref{subsection 2.1}
that $\mathcal B(\mathbf C)=\det(\mathbf C)\mathbf C^{-T}\mathbf L$ and multiplication of $\mathcal T$ with $\mathbf F^{T}$ in the third mode
yields a tensor of which every frontal slice is a rank-1 matrix, proportional to $\mathbf a_r\mathbf b_r^T$ for some $r\in\{1,\dots,R\}$.
On the other hand, if $\mathbf C$ is rectangular ($K<R$), then multiplication with $\mathbf F^T$ yields a tensor of which all slices are rank-$(R-K+1)$ matrices, generated by
$R-K+1$ rank-1 matrices $\mathbf a_r\mathbf b_r^T$. If we choose slices that have all but one rank-1 matrix in common, these form a tensor
that is as in Theorem \ref{theoremoralggeneigdec}
and of which the CPD yields $R-K+2$ columns of $\mathbf A$ and $\mathbf B$.
The result is formalized in the following lemma. The second statement implies that we do not have to compute the CPD to verify whether a slice combination is suitable.
\begin{lemma}\label{lemma:4.3}
Let   $\mathcal T=[\mathbf A,\mathbf B,\mathbf C]_R$,
the matrix  $\mathbf F\in\mathbb R^{K\times C^{K-1}_R}$ coincide with $\mathcal B(\mathbf C)$ up to column permutation and  scaling,
$\mathcal V=[\mathbf A,\mathbf B,\mathbf F^T\mathbf C]_R$,
$k_{\mathbf C}=K$, $m:= R-K+2$, and
$\mathcal C_m(\mathbf A)\odot \mathcal C_m(\mathbf B)$  have full column rank. Let also  $\mathbf y_1,\dots,\mathbf y_{C^{K-1}_R}$ denote the
 columns of $\mathbf C^T\mathbf F$ and
$\mathbf V_1,\dots,\mathbf V_{C^{K-1}_R}$ denote the frontal slices of $\mathcal V$.
Then the following statements are equivalent.
\begin{itemize}
\item[\textup{(i)}] The matrix $[\mathbf y_i\ \mathbf y_j]$ has exactly $K-2$ zero rows.
\item[\textup{(ii)}] The matrices $[\mathbf V_i\ \mathbf V_j]$ and  $[\mathbf V_i^T\ \mathbf V_j^T]$ have rank $m$.
\item[\textup{(iii)}] The  tensor $\mathcal V_{ij}$ formed by the  frontal slices $\mathbf V_i$ and $\mathbf V_j$ has rank $m$.
\end{itemize}
The CPD
$\mathcal V_{ij}=[\left[\begin{matrix}\mathbf a_{p_1}&\dots&\mathbf a_{p_m}\end{matrix}\right],\left[\begin{matrix}\mathbf b_{p_1}&\dots&\mathbf b_{p_m}\end{matrix}\right], \widehat{\mathbf C}]_m$ can be found algebraically by Theorem \ref{theoremoralggeneigdec} \textup{(ii)} and the indices $p_1,\dots,p_m$ are uniquely defined by
the pair $(i,j)$.
\end{lemma}
\begin{proof}
Since the proof is technical, it is given in the supplementary materials.
\end{proof}

To summarize, we first
find a matrix $\mathbf F\in\mathbb R^{K\times C^{K-1}_R}$  that
coincides with $\mathcal B(\mathbf C)$ up to column permutation and  scaling.
 We construct a tensor $\mathcal V$ from the tensor $\mathcal T$ and the matrix $\mathbf F$ as in Lemma  \ref{lemma:4.3} and choose
 slice combinations for which
   $[\mathbf V_i\ \mathbf V_j]$ and  $[\mathbf V_i^T\ \mathbf V_j^T]$ have rank $m$. For each such slice combination
 the CPD of the  corresponding tensor
yields $m$ columns of $\mathbf A$ and $\mathbf B$.
In this way we obtain   all columns of $\mathbf A$ and $\mathbf B$.
Overall, there exist exactly $C^m_RC^2_m$ pairs $(i,j)$ such that \textup{(i)}--\textup{(iii)} hold.
The amount of work  can be reduced by finding enough, instead of all, tensors $\mathcal V_{ij}$ that yield columns of $\mathbf A$ and $\mathbf B$.
The matrix $\mathbf C$ is finally obtained by  $\mathbf C=\left((\mathbf A\odot\mathbf B)^\dagger\Matr{\mathcal T}\right)^T$.
 \subsection{Discussion of working conditions}\label{Section:5}
One may wonder what happens if the conditions in Theorem \ref{theorem: main} are not satisfied, or, the other way around, under which circumstances the algorithms will fail. It turns out that, at least if the tensor rank is known, the crucial steps are steps 2 and 5 of Phase 1. If these do not pose problems, then the overall algorithms will work. Step 2 poses a problem when  
$\dim\ker(\TTm{m})\geq C^{K-1}_{R}$. This indicates  that $k_{\mathbf C}<K-1$ and/or that 
$\mathcal C_m(\mathbf A)\odot \mathcal C_m(\mathbf B)$  does not have full column rank. If step 2 does not pose a problem, but step 5 does, then
$k_{\mathbf C}=K-1$. This is formalized in the following lemma.
  \begin{lemma}
Let $\mathcal T=[\mathbf A,\mathbf B,\mathbf C]_R$ be a CPD of $\mathcal T$ and $m=R-K+2$. Let
the matrix $\TTm{m}$ be defined by Definition \ref{definition:detecting matrix}, and let the tensor $\mathcal W$ be constructed in
steps 3--4 of Phase 1.
\begin{itemize}
\item[\textup{(i)}] If  $\dim\ker(\TTm{m})=C^{K-1}_{R}$, then
$k_{\mathbf C}\geq K-1$ and the matrix $\mathcal C_m(\mathbf A)\odot \mathcal C_m(\mathbf B)$  has full column rank;
\item[\textup{(ii)}] If additionally,  $r_{\mathcal W}=C^{K-1}_R$ and 
$\mathcal W=[\mathbf F,\mathbf F_2,\mathbf F_3]_{C^{K-1}_R}$, where $k_{\mathbf F}\geq 2$ and the matrices
$\mathbf F_2$ and $\mathbf F_3$ have full column rank, then $k_{\mathbf C}=K$.
\end{itemize}
\end{lemma}
\begin{proof}
\textup{(i)} The proof is by contradiction. 
Assume that $k_{\mathbf C}\leq K-2$. Then $\mathbf C$ has
$K-1$ columns that are linearly dependent. Without loss of generality we can assume that these columns are $\mathbf c_m,\dots,\mathbf c_R$. Then the columns
$\{\pi_s(\mathbf c_1\otimes\dots\otimes \mathbf c_{m-1}\otimes\mathbf c_k)\}_{k=m}^R$  of the matrix $\mathcal R_m(\mathbf C)$ are also linearly dependent.
Hence, by Lemma \ref{lemma:2.22} \textup{(ii)}, the $C^m_R\times C^m_{K+m-1}$ matrix $\mathcal Q_m(\mathbf C)^T$ has linearly dependent rows, which implies that
$\dim \left( \ker \left(\mathcal{Q}_m(\mathbf C)^T\right)\right)\geq C^m_{K+m-1}-C^m_R+1=C^{K-1}_R+1$. On the other hand, by
Lemma \ref{lemma:3.9} \textup{(ii)}, $\dim(\ker (\TTm{m}))\geq\ \dim \left( \ker \left(\mathcal{Q}_m(\mathbf C)^T\right)\right)$,
 which is a contradiction with $\dim(\ker (\TTm{m}))=C^{K-1}_R$.

 We have proved that $k_{\mathbf C}\geq K-1$. 
By Proposition \ref{proposition:2.23} \textup{(i)}, the matrix $\mathcal R_m(\mathbf C)$ has full column rank.
Since by Lemma \ref{lemma:2.22} \textup{(i)}, $\mathcal R_m(\mathbf C)=\mathbf H^T\mathcal Q_m(\mathbf C)$, it follows that the matrix $\mathcal Q_m(\mathbf C)$ also has full column rank. Hence,
$\dim \left( \ker \left(\mathcal{Q}_m(\mathbf C)^T\right)\right)=C^m_{K+m-1}-C^m_R=C^{K-1}_R$. Hence, by Lemma
\ref{lemma:3.9} \textup{(ii)},
$$
\dim(\ker (\TTm{m}))=
\dim(\ker(\mathcal C_m(\mathbf A)\odot \mathcal C_m(\mathbf B)))+\dim \left( \ker \left(\mathcal{Q}_m(\mathbf C)^T\right)\right).
$$
Since, by assumption, $\dim\ker(\TTm{m})=C^{K-1}_{R}$, 
the matrix $\mathcal C_m(\mathbf A)\odot \mathcal C_m(\mathbf B)$ has full column rank.

\textup{(ii)} It is clear that  the columns of
 the matrix $\mathbf F\odot \mathbf F_2$ form  a basis of  $\textup{range}(\mathbf W)$. 
 By construction of $\mathbf W$,  $\textup{range}(\mathbf W)= \ker( \RTm{m}\upharpoonright_{\textup{range}(\pi_S)})$, and by
 Lemma \ref{lemma:3.9} \textup{(i)}, $\ker( \RTm{m}\upharpoonright_{\textup{range}(\pi_S)})=\ker\left(\mathcal{R}_m(\mathbf C)^T\upharpoonright_{\textup{range}(\pi_S)}\right)$. Hence, by Corollary \ref{corollary:2.18},
 $k_{\mathbf C}=K$.
\end{proof}
	\subsection{Theorem \ref{corrolary: main1}}
It remains to prove  Theorem \ref{corrolary: main1} (ii).
In the proof we construct a new tensor $\bar{\mathcal T}$ that has the same first two factor matrices as
$\mathcal T$ and the CPD of which  can be found by Algorithm \ref{alg:Main} or Algorithm \ref{alg:Main2}. 
Although, by construction of $\bar{\mathcal T}$, its frontal slices   are random linear combinations of
the frontal slices of $\mathcal T$, we still call  the overall procedure ``algebraic'' because
the proof of Theorem \ref{theoremoralggeneigdec} is also based on the same random slice mixture idea
(see \cite{Leurgans1993} and  references therein).

{\em Proof.}
Let the matrix $\mathbf C$ have $K$ rows, $\mathbf X$ be a $k_{\mathbf C}\times K$ matrix,  and
$\bar{\mathcal T}:=[\mathbf A,\mathbf B,\mathbf X\mathbf C]_R$. Then $\mathbf X\mathbf C\in \mathbb R^{k_{\mathbf C}\times R}$ and by  \eqref{eqT_V},
$
\Matr{\mathcal T}\mathbf X^T:=
(\mathbf A\odot \mathbf B)\mathbf C^T\mathbf X^T=
(\mathbf A\odot \mathbf B)(\mathbf X\mathbf C)^T=
\Matr{\bar{\mathcal T}}.
$
Thus, the multiplication of the third factor matrix of $\mathcal T$ by  $\mathbf X$ from the left is
equivalent to the multiplication of the matrix unfolding $\Matr{\mathcal T}$ by $\mathbf X^T$ from the right.

\textup{(i)}\ Assume that $\mathbf X$ is such that  $r_{\mathbf X\mathbf C}=k_{\mathbf X\mathbf C}=k_{\mathbf C}$. Then by Theorem \ref{theorem: main},
the CPD of $\bar{\mathcal T}$ is unique and can be found algebraically.
In particular, the matrix $\mathbf A\odot\mathbf B$ has full column rank and can be found up to column permutation and  scaling.
Hence, $\mathbf C=\left((\mathbf A\odot\mathbf B)^\dagger\Matr{\mathcal T}\right)^T$, and the proof is completed.

\textup{(ii)}\ It remains to present a construction of the matrix $\mathbf X$ such that
$k_{\mathbf X\mathbf C}=k_{\mathbf C}$. It is clear that  $k_{\mathbf X\mathbf C}\leq k_{\mathbf C}$. We claim that
$k_{\mathbf X\mathbf C}= k_{\mathbf C}$ for  generic $\mathbf X$. Namely,
\begin{equation}
\mu\{\vect{\mathbf X}:\ \mathbf X\in\mathbb R^{k_{\mathbf C}\times K},\ \ k_{\mathbf X\mathbf C}<k_{\mathbf C}   \}=0,\label{muzero}
\end{equation}
where  $\mu$ denotes the Lebesgue measure on $\mathbb R^{k_{\mathbf C}K}$.
It is well known that the zero set of a nonzero polynomial has Lebesgue measure  zero. Hence, for a
 nonzero vector $\mathbf f\in\mathbb R^{C^{k_{\mathbf C}}_K}$ we obtain
\begin{equation}
\mu\{\vect{\mathbf X}:\ \mathbf X\in\mathbb R^{k_{\mathbf C}\times K},\ \mathcal C_{k_{\mathbf C}}(\mathbf X)\mathbf f=0   \}=0.
\label{mu_f}
\end{equation}
From  Lemma \ref{compoundprop1} (1)
it follows that the matrix $\mathcal C_{k_{\mathbf C}}(\mathbf C)$ has all columns nonzero. By Lemma \ref{LemmaCompound}, $k_{\mathbf X\mathbf C}<k_{\mathbf C}$  if and only if the vector
$\mathcal C_{k_{\mathbf C}}(\mathbf X)\mathcal C_{k_{\mathbf C}}(\mathbf C)=\mathcal C_{k_{\mathbf C}}(\mathbf X\mathbf C)$ has a zero entry.
Hence, by \eqref{mu_f},
\begin{equation}
\begin{split}
&\{\vect{\mathbf X}:\ \mathbf X\in\mathbb R^{k_{\mathbf C}\times K},\ \ k_{\mathbf X\mathbf C}<k_{\mathbf C}   \}=\\
&\{\vect{\mathbf X}:\ \mathbf X\in\mathbb R^{k_{\mathbf C}\times K},\ \ \mathcal C_{k_{\mathbf C}}(\mathbf X)\mathcal C_{k_{\mathbf C}}(\mathbf C) \text{ has a zero entry}   \}=\\
&\bigcup\limits_{\mathbf f \text{ is a column of }\mathcal C_{k_{\mathbf C}}(\mathbf C)}
\{\vect{\mathbf X}:\ \mathbf X\in\mathbb R^{k_{\mathbf C}\times K},\ \ \mathcal C_{k_{\mathbf C}}(\mathbf X)\mathbf f=0\}.\label{mu_fmany}
\end{split}
\end{equation}
Now \eqref{muzero} follows from
\eqref{mu_fmany} and \eqref{mu_f}.\qquad\endproof

Let the conditions of Theorem \ref{corrolary: main1},
Corollary \ref{likegenKruskal} or Corollary \ref{corrolary: main2} hold.
The following procedure for computing the CPD
follows from the proof of Theorem \ref{corrolary: main1} (ii).
First, we generate  a random  $k_{\mathbf C}\times K$ matrix $\mathbf X$ and set
$\bar{\mathbf T}=\Matr{\mathcal T}\mathbf X^T$.
 Then  the $I\times J\times k_{\mathbf C}$ tensor
 $\bar{\mathcal T}:=\UnMatr{\bar{\mathbf T},I,J}$
 satisfies the conditions of Theorem \ref{theorem: main} (ii) with $K$ replaced by $k_{\mathbf C}$. Hence, the CPD  $\bar{\mathcal T}=[\mathbf A, \mathbf B,\mathbf X\mathbf C]_R$  can be found by Algorithm \ref{alg:Main} or Algorithm \ref{alg:Main2}. Finally,
 the matrix $\mathbf C$ is obtained by  $\mathbf C=\left((\mathbf A\odot\mathbf B)^\dagger\Matr{\mathcal T}\right)^T$.
 \subsection{Examples}
\begin{example}
Let  $\mathcal T=[\mathbf A,\mathbf B,\mathbf C]_5$ with
$$
\mathbf A=
\left[\begin{matrix}
        1 &1& 0& 0& 0\\
        1 &0& 1& 0& 0\\
        1 &0 &0 &1 &0\\
        0 &0 &0 &0 &1
        \end{matrix}
        \right],\quad
\mathbf B=
\left[\begin{matrix}
1& 0& 0& 0& 1\\
        1& 0& 0& 1& 0\\
        1& 0& 1& 0& 0\\
        0& 1& 0& 0& 0
\end{matrix}
        \right],\quad
\mathbf C=\left[\begin{matrix}
        1& 1& 0& 0& 0\\
        1& 0& 2& 0& 0\\
        1& 0& 0& 3& 0\\
        1& 0& 0& 0& 1
        \end{matrix}
        \right].
$$
Since condition  \eqref{eqtwomatrK2} does not hold, the rank and uniqueness of the CPD do not follow from
Kruskal's Theorem \ref{theoremKruskalnew1}.
One can easily check that the conditions of Theorem \ref{theorem: main} hold for
$m=5-4+2=3$.
Hence, the factor matrices of $\mathcal T$ can be found
by Algorithms \ref{alg:Main} and \ref{alg:Main2}.

{\bf Phase 1 of Algorithms \ref{alg:Main} and \ref{alg:Main2}.}
The  frontal slices of $\mathcal T$ are
$$
\mathbf T_1=\left[\begin{matrix}1&1&1&1\\1&1&1&0\\1&1&1&0 \\0 &0 &0 &0 \end{matrix}\right],\quad
\mathbf T_2=\left[\begin{matrix}1&1&1&0\\1&1&1&0\\1&4&1&0 \\0 &0 &0 &0 \end{matrix}\right],\quad
\mathbf T_3=\mathbf T_1^T,\quad
\mathbf T_4=\left[\begin{matrix}1&1&1&0\\1&1&3&0\\1&1&1&0 \\0 &0 &0 &0
\end{matrix}\right].
$$
We construct  the $C^3_4C^3_4$-by-$C^3_{6}$ (or $16$-by-$20$) matrix $\TTm{3}$ by Definition \ref{definition:detecting matrix}. For instance,
the $(1,2,3)$-rd (or the $6$-th) column   of  $\TTm{3}$ is equal to
$\vect{\rdm{2}(\mathbf T_1,\mathbf T_2,\mathbf T_3)}$, where $\rdm{2}(\mathbf T_1,\mathbf T_2,\mathbf T_3)$ is computed by
\eqref{eq:example123} and equals
$$
\rdm{2}(\mathbf T_1,\mathbf T_2,\mathbf T_3)=
-\left[
\begin{array}{rrrr}
     0&    -3&     0&     3\\
     0&     1&     1&     0\\
     3&     4&     1&     0\\
     3&     0&     0&     0
\end{array}
\right].
$$
The full  matrix $\TTm{3}$ is given in
the supplementary materials.
It can be  checked that $\ker (\TTm{3})=\textup{range}(\overline{\mathbf W})$, where
$$
\overline{\mathbf W}=[\begin{matrix}\mathbf e_1^{20}& \mathbf e_{11}^{20}& \mathbf e_{17}^{20}& \mathbf e_{20}^{20}& \mathbf e_{2,-5}^{20}& \mathbf e_{4,-10}^{20}& \mathbf e_{3,-8}^{20}& \mathbf e_{13,-16}^{20}& \mathbf e_{12,-14}^{20}& \mathbf e_{18,-19}^{20}\end{matrix}]
$$
and $\mathbf e_{i,-j}^{20}:=\mathbf e_{i}^{20}-\mathbf e_{j}^{20}$.  Let $\mathbf G$ be the $64\times 20$ matrix   defined by \eqref{eq:2.11}.
We denote by $\mathcal W$ the $4\times 16\times 10$ tensor such that $\Matr{\mathcal W}=\mathbf G\overline{\mathbf W}$.
We find algebraically the CPD $\mathcal W=[\mathbf F,\mathbf F_2,\mathbf F_3]_{10}$ with
\begin{gather*}
\mathbf F=
\left[
\begin{array}{rrrrrrrrrr}
1&-1&0&-1&0&0&0&-1&0&0\\
0&1&0&0&-1&1&0&0&0&-1\\
0&0&0&1&0&-1&-1&0&-1&0\\
1&0&1&0&0&0&0&0&1&1
\end{array}
\right],
\end{gather*}
$\mathbf F_2=\mathbf F\odot\mathbf F$ and some nonsingular matrix $\mathbf F_3$.
In the sequel we will use only the fact that $\mathbf F$ coincides with $\mathcal B(\mathbf C)$ up to column permutation and scaling.

{\bf Phase 2 and 3 of Algorithm \ref{alg:Main}.}
There are
$210$ $4\times 6$  submatrices of $\mathbf F$.
In Phase 2 of Algorithm 1 we  pick the five submatrices  that have rank  $3$. One can easily see that these submatrices are $[\mathbf f_1\ \mathbf f_2\ \mathbf f_3\ \mathbf f_5\ \mathbf f_8\ \mathbf f_{10}]$, $[\mathbf f_1\ \mathbf f_2\ \mathbf f_4\ \mathbf f_6\ \mathbf f_9\ \mathbf f_{10}]$, $[\mathbf f_1\ \mathbf f_3\ \mathbf f_4\ \mathbf f_7\ \mathbf f_8\ \mathbf f_{9}]$, $[\ \mathbf f_2\ \mathbf f_4\ \mathbf f_5\ \mathbf f_6\ \mathbf f_7\ \mathbf f_8]$, $[\mathbf f_3\ \mathbf f_5\ \mathbf f_6\ \mathbf f_7\ \mathbf f_9\ \mathbf f_{10}]$. Their left kernels have dimension $1$ and are spanned by the norm one vectors
$\widehat{\mathbf c}_1=[0\ 0\ 1\ 0]^T$, $\widehat{\mathbf c}_2=[0.5\ 0.5\ 0.5\ 0.5]^T$, $\widehat{\mathbf c}_3=[0\ 1\ 0\ 0]^T$, $\widehat{\mathbf c}_4=[0\ 0\ 0\ 1]^T$, $\widehat{\mathbf c}_5=[1\ 0\ 0\ 0]^T$, respectively. The matrix formed by these vectors coincides with the matrix  $\mathbf C$ up to column permutation and scaling.

Let us demonstrate how  Phase 3 of Algorithm \ref{alg:Main} works. One can easily see that the vectors
$\mathbf z_1=[-1\ 0\ 1\ 0]^T$ and $\mathbf z_2=[2\ 0\ 0\ 0]^T$ coincide with the first two columns of the matrix
$[\widehat{\mathbf c}_1\ \widehat{\mathbf c}_2\ \widehat{\mathbf c}_3\ \widehat{\mathbf c}_4]^{-T}$.
We have
$$
\Matr{\widetilde{\mathcal T}}=\Matr{\mathcal T}[\mathbf z_1\ \mathbf z_2]=(\mathbf A\odot\mathbf B)\mathbf C^T
[\mathbf z_1\ \mathbf z_2]=
(\mathbf A\odot\mathbf B)\left[\begin{array}{rrrrr}
0&  -1&  0&  3&  0\\
2&   2&  0&  0&  0\end{array}\right]^T.
$$
or $\widetilde{\mathcal T}=\left[[\mathbf a_1\ \mathbf a_2\ \mathbf a_4], [\mathbf b_1\ \mathbf b_2\ \mathbf b_4], \left[\begin{array}{rrr}0&-1&3\\2&2&0\end{array}\right]\right]_3$.
Thus, computing algebraically the CPD of $\widetilde{\mathcal T}$ we find the vectors $\mathbf a_1\otimes\mathbf b_1$, $\mathbf a_2\otimes\mathbf b_2$, and $\mathbf a_4\otimes\mathbf b_4$. The vectors $\mathbf a_3\otimes\mathbf b_3$ and  $\mathbf a_5\otimes\mathbf b_5$ are found by
\begin{equation*}
[\mathbf a_3\otimes\mathbf b_3\ \ \mathbf a_5\otimes\mathbf b_5]=\left(
\Matr{\mathcal T}-[\mathbf a_1\otimes\mathbf b_1\ \ \mathbf a_2\otimes\mathbf b_2\ \ \mathbf a_4\otimes\mathbf b_4][\mathbf c_1\ \ \mathbf c_2\ \ \mathbf c_4]^T\right)[\mathbf c_3\ \ \mathbf c_5]^{\dagger,T}.
\end{equation*}
{\bf Phase 2 and 3 of Algorithm \ref{alg:Main2}.}
We construct the $4\times 4\times 10$ tensor $\mathcal V$ with matrix unfolding $\Matr{\mathcal V}=\Matr{\mathcal T}\mathbf F$.
Let $\mathbf V_1,\dots,\mathbf V_{10}$ denote the  frontal slices of $\mathcal V$ and let $\mathcal V_{ij}$ denote the $4\times 4\times 2$ tensor with frontal slices
$\mathbf V_i$ and $\mathbf V_j$. We construct the set
\begin{equation*}
\begin{split}
\mathcal J:=& \{(i,j): \text{ the matrices }[\mathbf V_i\ \mathbf V_j] \text{ and  }[\mathbf V_i^T\ \mathbf V_j^T]\text{ have rank }3,\ 1\leq i<j\leq 10\}=\\
&\left\{ (1,2),\ (1,3),\ (1,4),\ (1,8),\ (1,9),\ (1,10),\ (2,4),\ (2,5),\ (2,6),\ (2,8),\ (2,10),\right.\\
&\ \  (3,5),\ (3,7),\ (3,8),\ (3,9),\ (3,10),\ (4,6),\ (4,7),\ (4,8),\ (4,9),\\
&\ \ \left. (5,6),\ (5,7),\ (5,8),\ (5,10),\ (6,7),\ (6,9),\ (6,10),\ (7,8),\ (7,9),\ (9,10)\right\}.
\end{split}
\end{equation*}
For $(i,j)\in \mathcal J$, $\mathcal V_{ij}$ has rank $3$ and the CPD can be computed algebraically.
For instance,
$
\Matr{\mathcal V_{12}}=\Matr{\mathcal T}[\mathbf f_1\ \mathbf f_2]=(\mathbf A\odot\mathbf B)\mathbf C^T[\mathbf f_1\ \mathbf f_2]
$. Since
$
\mathbf C^T[\mathbf f_1\ \mathbf f_2]=$
$
\left[\begin{array}{rrrrr} 2&1&0&0&0\\0&-1&2&0&0\end{array}\right]^T
$,
 we have $\mathcal V_{12}=\left[[\mathbf a_1\ \mathbf a_2\ \mathbf a_3], [\mathbf b_1\ \mathbf b_2\ \mathbf b_3], \left[\begin{array}{rrr}2&1&0\\0&-1&2\end{array}\right]\right]_3$.
In this way for each pair $(i,j)\in \mathcal J$  we estimate up to column scaling three columns of $\mathbf A$ and the corresponding columns of $\mathbf B$.
If we store all the estimates of columns of $\mathbf A$ and $\mathbf B$ in  $4\times 90$ matrices $\widetilde{\mathbf A}$ and $\widetilde{\mathbf B}$
then $\widetilde{\mathbf A}\odot\widetilde{\mathbf B}$  will contain $5$ clusters of $18$  collinear columns. Taking the cluster centers we get a matrix  $\mathbf Z$ which coincides with $\mathbf A\odot\mathbf B$ up to column scaling and permutation. Finally, the matrix
$\left(\mathbf Z^\dagger\Matr{\mathcal T}\right)^T$ coincides with $\mathbf C$ up to column scaling and the same permutation.
\end{example}
\begin{example}\label{remark
last}
It was shown in \cite{PartII} that the conditions of Theorem
\ref{theorem: main} hold for a generic  $6\times 6\times 7$ tensor of rank $9$.
This case is beyond Kruskal's bound.
Let $\mathbf F$ be the $7\times 84$ matrix produced by  Phase 1 of Algorithms \ref{alg:Main} and \ref{alg:Main2}.
Each column of the third factor matrix of the tensor is orthogonal to exactly $42$ columns of the matrix $\mathbf F$.
Since $C^{42}_{84}$ is of order $10^{24}$,  Phase 2  as presented in Algorithm \ref{alg:Main}  is computationally infeasible.
On the other hand, in  Phase 2 of Algorithm \ref{alg:Main2} we check the rank of  $2C^2_{84}=3486$ matrices of the size $6\times 12$ each.
 Then we have to  find algebraically the CPD of  $C^4_9C^2_4=756$ rank-$4$ tensors with dimensions $6\times 6\times 2$, which is
 equivalent with the computation of the GEVD of the associated matrix pencils. Moreover, one may further limit the amount of work by only determining subsets of $\mathcal J$.
We implemented Algorithm \ref{alg:Main2} in
MATLAB 2008a and we did experiments on a computer
with Intel(R) Core(TM) T9600 Duo 2.80GHz CPU and 4GB memory running
Windows Vista. The simulations demonstrate that
with a suboptimal implementation, it takes less than
$9$ seconds to compute the CPD of a generic  $6\times 6\times 7$ tensor of rank $9$.
\end{example}
\section{Conclusion}
We have proposed two  algorithms to compute CPD.
Both algorithms are algebraic  in the sense that they rely only on standard linear algebra and reduce the problem to the computation of GEVD.
The reduction exploits properties of
(polarized) compound matrices and permanents.
 The derivation spans the possibilities from \cite{Kruskal1977} to \cite{DeLathauwer2006,JiangSid2004} and covers cases beyond Kruskal's bound.

In this paper we have limited ourselves to exact CPD. In applications, CPD most often only approximates the given (noisy) tensor. A first observation is that the ``exact result'' could be used to initialize iterative algorithms for problem \eqref{eq:optproblem}.
We also note that \eqref{eq:newCPD} may be interpreted as the CPD of a partially symmetric tensor of order $m+1$ of which the first $m$ factor matrices are equal and parameterized by $\mathbf C$. This is a structure that can be handled by current algorithms in Tensorlab \cite{Tensorlab}. These algorithms are optimization-based and are not formally guaranteed to find the solution. However, they show excellent performance in practice. So far, we have computed
$ \ker(\RTm{m}\upharpoonright_{\textup{range}(\pi_S)})$
and then we have fitted the CPD structure to the result. Numerically, we could go a step further and take the Khatri-Rao structure into account in the computation of the kernel itself, with the kernel vectors parameterized by $\mathbf C$ and $\mathbf M$. One may also investigate whether the Khatri-Rao structure and the structure of $\RTm{m}$ may be exploited to avoid the computation of the mixed discriminants, so that one obtains an algorithm that works directly on $\mathcal T$. Since numerical aspects lead to a different type of study, we choose to defer them to an other paper.


\begin{thebibliography}{10}

\bibitem{Aleksandrov1938}
{\sc A.~D. Aleksandrov}, {\em Zur {T}heorie der gemischten {V}olumina von
  konvexen {K}\"{o}rpern. \textup{IV}. {D}ie gemischten {D}iskriminanten und
  die gemischten {V}olumina ({R}ussian) ({G}erman summary)}, Mat. Sb., 3 (45),
  no. 2 (1938), pp.~227--251.

\bibitem{Bapat1989}
{\sc R.B. Bapat}, {\em Mixed discriminants of positive semidefinite matrices},
  Linear Algebra Appl., 126 (1989), pp.~107--124.

\bibitem{1970_Carroll_Chang}
{\sc J.~Carroll and J.-J. Chang}, {\em {A}nalysis of individual differences in
  multidimensional scaling via an {N}-way generalization of
  ``{E}ckart-{Y}oung'' decomposition}, Psychometrika, 35 (1970), pp.~283--319.

\bibitem{LievenCichocki2013}
{\sc A.~Cichocki, D.~Mandic, C.~Caiafa, A-H. Phan, G.~Zhou, Q.~Zhao, and
  L.~De~Lathauwer}, {\em {T}ensor {D}ecompositions for {S}ignal {P}rocessing
  {A}pplications. {F}rom {T}wo-way to {M}ultiway {C}omponent {A}nalysis},
  ESAT-STADIUS Internal Report, 13-235, Leuven, Belgium: Department of
  Electrical Engineering (ESAT), KU Leuven,  (2013).

\bibitem{2009Comonetall}
{\sc P.~Comon, X.~Luciani, and A.~L.~F. de~Almeida}, {\em Tensor
  decompositions, alternating least squares and other tales}, J. Chemometrics,
  23 (2009), pp.~393--405.

\bibitem{DeLathauwer2006}
{\sc L.~De~Lathauwer}, {\em A {L}ink {B}etween the {C}anonical {D}ecomposition
  in {M}ultilinear {A}lgebra and {S}imultaneous {M}atrix {D}iagonalization},
  SIAM J. Matrix Anal. Appl., 28 (2006), pp.~642--666.

\bibitem{Lieven_ISPA}
\leavevmode\vrule height 2pt depth -1.6pt width 23pt, {\em A short introduction
  to tensor-based methods for factor analysis and blind source separation}, in
  ISPA 2011: Proceedings of the 7th International Symposium on Image and Signal
  Processing and Analysis,  (2011), pp.~558--563.

\bibitem{PartI}
{\sc I.~Domanov and L.~De~Lathauwer}, {\em {O}n the {U}niqueness of the
  {C}anonical {P}olyadic {D}ecomposition of {T}hird-{O}rder {T}ensors ---
  {P}art {I}: {B}asic {R}esults and {U}niqueness of {O}ne {F}actor {M}atrix},
  SIAM J. Matrix Anal. Appl., 34 (2013), pp.~855--875.

\bibitem{PartII}
\leavevmode\vrule height 2pt depth -1.6pt width 23pt, {\em {O}n the
  {U}niqueness of the {C}anonical {P}olyadic {D}ecomposition of {T}hird-{O}rder
  {T}ensors--- {P}art {II}: {O}verall {U}niqueness}, SIAM J. Matrix Anal.
  Appl., 34 (2013), pp.~876--903.

\bibitem{Egorychev1981}
{\sc G.~P. Egorychev}, {\em {P}roof of the van der {W}aerden conjecture for
  permanents}, Siberian Math. J., 22 (1981), pp.~854–--859.

\bibitem{Harshman1970}
{\sc R.~A. Harshman}, {\em Foundations of the {PARAFAC} procedure: {M}odels and
  conditions for an "explanatory" multi-modal factor analysis}, UCLA Working
  Papers in Phonetics, 16 (1970), pp.~1--84.

\bibitem{Harshman1972}
\leavevmode\vrule height 2pt depth -1.6pt width 23pt, {\em Determination and
  {P}roof of {M}inimum {U}niqueness {C}onditions for{ PARAFAC}1}, UCLA Working
  Papers in Phonetics, 22 (1972), pp.~111--117.

\bibitem{1994HarshmanLundy}
{\sc R.~A. Harshman and M.~E. Lundy}, {\em {P}arafac: {P}arallel factor
  analysis}, Comput. Stat. Data Anal.,  (1994), pp.~39--72.

\bibitem{Hitchcock}
{\sc F.~L. Hitchcock}, {\em The expression of a tensor or a polyadic as a sum
  of products}, J. Math. Phys., 6 (1927), pp.~164--189.

\bibitem{HornJohnson}
{\sc R.~A. Horn and C.~R. Johnson}, {\em {M}atrix {A}nalysis}, Cambridge
  University Press, 1990.

\bibitem{JiangSid2004}
{\sc T.~Jiang and N.~D. Sidiropoulos}, {\em {K}ruskal's {P}ermutation {L}emma
  and the {I}dentification of {CANDECOMP}/{PARAFAC} and {B}ilinear {M}odels
  with {C}onstant {M}odulus {C}onstraints}, IEEE Trans. Signal Process., 52
  (2004), pp.~2625--2636.

\bibitem{KoldaReview}
{\sc T.~G. Kolda and B.~W. Bader}, {\em {T}ensor {D}ecompositions and
  {A}pplications}, SIAM Review, 51 (2009), pp.~455--500.

\bibitem{Kroonenberg2008}
{\sc P.~M Kroonenberg}, {\em {A}pplied {M}ultiway {D}ata {A}nalysis}, Hoboken,
  NJ: Wiley, 2008.

\bibitem{Kruskal1977}
{\sc J.~B. Kruskal}, {\em Three-way arrays: rank and uniqueness of trilinear
  decompositions, with application to arithmetic complexity and statistics},
  Linear Algebra Appl., 18 (1977), pp.~95--138.

\bibitem{Landsberg}
{\sc J.~M. Landsberg}, {\em {T}ensors: {G}eometry and {A}pplications}, AMS,
  Providence, Rhode Island, 2012.

\bibitem{Leurgans1993}
{\sc S.~E. Leurgans, R.~T. Ross, and R.~B. Abel}, {\em A decomposition for
  three-way arrays}, SIAM J. Matrix Anal. Appl., 14 (1993), pp.~1064--1083.

\bibitem{determinantfunctions}
{\sc J.~S. Lomont and M.~S. Cheena}, {\em A multilinearity property of
  determinant functions}, Linear and Multilinear Algebra, 14 (1983),
  pp.~199--223.

\bibitem{Marcus1964}
{\sc M.~Marcus}, {\em The {H}adamard {T}heorem for {P}ermanents}, Proc. Amer.
  Math. Soc., 15 (1964), pp.~pp. 967--973.

\bibitem{marcus1973finite}
\leavevmode\vrule height 2pt depth -1.6pt width 23pt, {\em Finite dimensional
  multilinear algebra}, no.~pt. 1 in Pure and applied mathematics, M. Dekker,
  1973.

\bibitem{MarcusMinc1961}
{\sc M.~Marcus and H.~Minc}, {\em On the relation between the determinant and
  the permanent}, Illinois J. Math., 5 (1961), pp.~376--381.

\bibitem{1988Topographic}
{\sc J.~M$\ddot{\text{o}}$cks}, {\em Topographic components model for
  event-related potentials and some biophysical considerations}, IEEE Trans.
  Biomed. Eng., 35 (1988), pp.~482--484.

\bibitem{Minc1978}
{\sc H.~Minc}, {\em Permanents}, Encyclopedia of Mathematics and its
  Applications, Addison-Wesley Publishing Company, Advanced Book Program, 1978.

\bibitem{Muir1882}
{\sc T.~Muir}, {\em A treatise on the theory of determinants: with graduated
  sets of exercises for use in colleges and schools}, Macmillan and Co., 1882.

\bibitem{DimLievenLL1}
{\sc D.~Nion and L.~De~Lathauwer}, {\em {A} {S}tudy of the {D}ecomposition of a
  {T}hird-{O}rder {T}ensor in {R}ank-$({L},{L},1)$ {T}erms}, ESAT-STADIUS
  Internal Report, 11-239, Leuven, Belgium: Department of Electrical
  Engineering (ESAT), KU Leuven,  (2011).

\bibitem{Oeding2013}
{\sc Luke Oeding and Giorgio Ottaviani}, {\em Eigenvectors of tensors and
  algorithms for {W}aring decomposition}, Journal of Symbolic Computation, 54
  (2013), pp.~9--35.

\bibitem{Sanchez1990}
{\sc E.~Sanchez and B.. Kowalski}, {\em Tensorial resolution: A direct
  trilinear decomposition}, J. Chemometrics, 4 (1990), pp.~29--45.

\bibitem{Sands1980}
{\sc R:~Sands and F.~Young}, {\em Component models for three-way data: An
  alternating least squares algorithm with optimal scaling features},
  Psychometrika, 45 (1980), pp.~39--67.

\bibitem{Schott2003}
{\sc James~R. Schott}, {\em Kronecker product permutation matrices and their
  application to moment matrices of the normal distribution}, J. Multivar.
  Anal., 87 (2003), pp.~177 -- 190.

\bibitem{smilde2004multi}
{\sc A.K. Smilde, R.~Bro, and P.~Geladi}, {\em Multi-way analysis with
  applications in the chemical sciences}, J. Wiley, 2004.

\bibitem{Tensorlab}
{\sc L.~Sorber, M.~Van~Barel, and L.~De~Lathauwer}, {\em Tensorlab v1.0},
  Available online, February 2013. URL:
  http://esat.kuleuven.be/sista/tensorlab/.

\bibitem{Sorber}
\leavevmode\vrule height 2pt depth -1.6pt width 23pt, {\em Optimization-based
  algorithms for tensor decompositions: canonical polyadic decomposition,
  decomposition in rank-(${L}_r$,${L}_r$,1) terms and a new generalization},
  SIAM J. Optim., 23 (2013), pp.~695–--720.

\bibitem{MikaelVDM}
{\sc M.~S{\o}rensen and L.~De~Lathauwer}, {\em {B}lind {S}ignal {S}eparation
  via {T}ensor {D}ecomposition with {V}andermonde {F}actor {P}art {I}:
  {C}anonical {P}olyadic {D}ecomposition}, IEEE Trans. Signal Process., 61
  (2013), pp.~5507--5519.

\bibitem{MikaelHO}
\leavevmode\vrule height 2pt depth -1.6pt width 23pt, {\em {N}ew {U}niqueness
  {C}onditions for the {C}anonical {P}olyadic {D}ecomposition of
  {H}igher-{O}rder {T}ensors}, ESAT-STADIUS Internal Report, 11-239, Leuven,
  Belgium: Department of Electrical Engineering (ESAT), KU Leuven,  (2013).

\bibitem{Mikaelorth}
{\sc M.~S{\o}rensen, L.~De~Lathauwer, P.~Comon, S.~Icart, and L.~Deneire}, {\em
  {C}anonical {P}olyadic {D}ecomposition with a {C}olumnwise {O}rthonormal
  {F}actor {M}atrix}, SIAM J. Matrix Anal. Appl., 33 (2012), pp.~1190--1213.

\bibitem{TenBerge2009}
{\sc J.~Ten~Berge and J.~Tendeiro}, {\em The link between sufficient conditions
  by {H}arshman and by {K}ruskal for uniqueness in {C}andecomp/{P}arafac}, J.
  Chemometrics, 23 (2009), pp.~321--323.

\end{thebibliography}
\end{document}


\pagestyle{myheadings}
\thispagestyle{plain}
\markboth{IGNAT DOMANOV AND LIEVEN DE LATHAUWER}{Supplementary materials}

\begin{center}\MakeUppercase{\bf Canonical Polyadic Decomposition of third-order tensors: reduction to generalized eigenvalue decomposition}
 \\ \qquad
 \\
  \MakeUppercase{Ignat Domanov
and
Lieven De Lathauwer}\\
\qquad\\
\MakeLowercase{Supplementary materials}
\end{center}

\section{Supplementary material related to Proposition 1.10}\qquad\\
Recall that the $K\times C^{K-1}_R$ matrix
$\mathcal B(\mathbf C)$ is defined by
\begin{equation}
\mathcal B(\mathbf C):=\mathbf L\mathcal C_{K-1}(\mathbf C),\label{def:B(C)}
\end{equation}
where
\begin{equation}
\mathbf L:=\left[\begin{array}{rrrc} 0     & 0     & \dots               &(-1)^{K-1}\\
                                     \vdots& \vdots& \iddots           &\vdots\\
                                0     & -1    & \dots               &0\\
                                1     & 0     & \dots                  &0
                 \end{array}\right].
                 \label{eq:matrixL}
\end{equation}
 We start with a trivial Lemma.
\begin{lemma}\label{lemma:Laplace}
Let $\mathbf x\in\mathbb R^K$, $\mathbf C\in\mathbb R^{K\times R}$, $k_{\mathbf C}\geq K-1$, and let $\mathcal B(\mathbf C)$ be defined by \eqref{def:B(C)}.
Then
\begin{itemize}
\item[\textup{(i)}]
\begin{equation*}
\mathbf x^T\mathcal B(\mathbf C)=
\left[\begin{matrix}
\det \left[\begin{matrix}\mathbf c_1&\dots& \mathbf c_{K-1}&\mathbf x\end{matrix}\right]
&\dots&
\det \left[\begin{matrix}\mathbf c_{R-K+2}&\dots& \mathbf c_{R}&\mathbf x\end{matrix}\right]
\end{matrix}\right];
\end{equation*}
\item[\textup{(ii)}]
$\mathcal B(\mathbf C)$ has no zero columns;
\item[\textup{(iii)}]
the  $(i_1,\dots,i_{K-1})$-th column of  $\mathcal B(\mathbf C)$  is orthogonal to
 $\textup{span}\{\mathbf c_{i_1},\dots,\mathbf c_{i_{K-1}}\}$ for $(i_1,\dots,i_{K-1})\in S^{K-1}_R$;
\item[\textup{(iv)}] if $k_{\mathbf C}=K$, then
$\mathcal B(\mathbf C)$ has no proportional columns, that is $k_{\mathcal B(\mathbf C)}\geq 2$.
\end{itemize}
\end{lemma}
\begin{proof}
\textup{(i)}\quad
From \eqref{def:B(C)} it follows that
the  $(i_1,\dots,i_{K-1})$-th column of  $\mathcal B(\mathbf C)$ is equal to
$\mathbf L\mathcal C_{K-1}\left[\begin{matrix}\mathbf c_{i_1}&\dots&\mathbf c_{i_{K-1}}\end{matrix}\right]$.
By  the Laplace expansion theorem
\begin{equation}\label{proofofiii}
\mathbf y^T\mathbf L\mathcal C_{K-1}\left[\begin{matrix}\mathbf c_{i_1}&\dots&\mathbf c_{i_{K-1}}\end{matrix}\right]=\det\left[\begin{matrix}\mathbf c_{i_1}&\dots&\mathbf c_{i_{K-1}}&\mathbf y\end{matrix}\right],\qquad \mathbf y\in\mathbb R^K.
\end{equation}
Now, the statement  \textup{(i)}  follows from \eqref{proofofiii} by setting  $\mathbf y=\mathbf x$.

\textup{(ii)}\quad
Since the vectors $\mathbf c_{i_1},\dots,\mathbf c_{i_{K-1}}$ are linearly independent in $\mathbb R^K$, it follows that there exists a vector $\mathbf y$
such that $\det\left[\begin{matrix}\mathbf c_{i_1}&\dots&\mathbf c_{i_{K-1}}&\mathbf y\end{matrix}\right]\ne 0$. Hence, by \eqref{proofofiii},
the  $(i_1,\dots,i_{K-1})$-th column of  $\mathcal B(\mathbf C)$  is nonzero.

\textup{(iii)} follows from \eqref{proofofiii} and the fact that $\det\left[\begin{matrix}\mathbf c_{i_1}&\dots&\mathbf c_{i_{K-1}}&\mathbf y\end{matrix}\right]=0$
if and only if $\mathbf y\in\textup{span}\{\mathbf c_{i_1},\dots,\mathbf c_{i_{K-1}}\}$.

\textup{(iv)}  Assume that the $(i_1,\dots,i_{K-1})$-th and the $(j_1,\dots,j_{K-1})$-th  column of $\mathcal B(\mathbf C)$ are proportional to a nonzero vector $\mathbf t$. Set  $E:=\textup{span}\{\mathbf c_{i_1},\dots,\mathbf c_{i_{K-1}},\mathbf c_{j_1},\dots,\mathbf c_{j_{K-1}}\}$.
Since $k_{\mathbf C}=K$, it follows that $E=\mathbb R^K$. On the other hand, by \textup{(iii)}, the nonzero vector $\mathbf t$ is orthogonal to $E$, which is a contradiction.
\end{proof}
\begin{lemma}\label{lemma:2.4}
Let $\mathbf C\in\mathbb R^{K\times R}$, $K\leq R$, and $k_{\mathbf C}=K$. Let also $\mathcal B(\mathbf C)$ be defined by \eqref{def:B(C)}. Then
\begin{itemize}
\item[\textup{(i)}] every column of the $R\times C_R^{K-1}$ matrix $\mathbf C^T\mathcal B(\mathbf C)$ has exactly $K-1$ zero entries. Namely, if
$\mathbf d=\left[\begin{matrix}d_1&\dots&d_R\end{matrix}\right]$ is the $(j_1,\dots,j_{K-1})$-th column of $\mathbf C^T\mathcal B(\mathbf C)$, then $d_{r}=0$ if and only if $r\in\{j_1,\dots,j_{K-1}\}$;
\item[\textup{(ii)}] every row of the $R\times C_R^{K-1}$ matrix $\mathbf C^T\mathcal B(\mathbf C)$ has exactly $C^{K-2}_{R-1}$ zero entries. Namely, if
$\widehat{\mathbf d}=\left[\begin{matrix}d_{(1,\dots,K-1)}&\dots&d_{(R-K+2,\dots,R)}\end{matrix}\right]$ is the $r$-th row of $\mathbf C^T\mathcal B(\mathbf C)$, then
$d_{(j_1,\dots,j_{K-1})}=0$ if and only if $r\in\{j_1,\dots,j_{K-1}\}$;
\item[\textup{(iii)}]   the matrix
$$
\mathcal B(\mathbf C)^{(l)}:= \underbrace{\mathcal B(\mathbf C)\odot\dots\odot \mathcal B(\mathbf C)}_{l}
$$
  has full column rank for  $l\geq R-K+1$.
\end{itemize}
\end{lemma}
\begin{proof}
\textup{(i)-(ii)}\quad By Lemma \ref{lemma:Laplace} \textup{(i)},
\begin{equation}
\begin{split}
&\mathbf C^T\mathcal B(\mathbf C)=\mathbf C^T\mathbf L\mathcal C_{K-1}(\mathbf C)=\\
&\left[\begin{matrix}
\det \left[\begin{matrix}\mathbf c_1&\dots& \mathbf c_{K-1}&\mathbf c_1\end{matrix}\right]
&\dots&
\det \left[\begin{matrix}\mathbf c_{R-K+2}&\dots& \mathbf c_{R}&\mathbf c_1\end{matrix}\right]\\
\vdots& \vdots &\vdots\\
\det \left[\begin{matrix}\mathbf c_1&\dots& \mathbf c_{K-1}&\mathbf c_R\end{matrix}\right]
&\dots&
\det \left[\begin{matrix}\mathbf c_{R-K+2}&\dots& \mathbf c_{R}&\mathbf c_R\end{matrix}\right]
\end{matrix}\right].
\end{split}\label{eq:2.8}
\end{equation}
Since $k_{\mathbf C}=K$, it follows that
\begin{equation}
\det \left[\begin{matrix}\mathbf c_{j_1}&\dots& \mathbf c_{j_{K-1}}&\mathbf c_r\end{matrix}\right]=0\ \Leftrightarrow
r\in\{j_1,\dots,j_{K-1}\}.\label{eq:2.9}
\end{equation}
The results now easily follow from \eqref{eq:2.8}-\eqref{eq:2.9}.

\textup{(iii)}\quad
It is sufficient to consider the case $l=R-K+1$. The result for $l\geq R-K+1$ then follows directly from the definition of the Khatri-Rao product.

Suppose that $\mathcal B(\mathbf C)^{(R-K+1)}\widehat{\mathbf t}=\vzero$ for
$\widehat{\mathbf t}=\left[\begin{matrix} t_{(1,\dots,K-1)}&\dots&  t_{(R-K+2,\dots,R)}\end{matrix}\right]\in\mathbb R^{C^{K-1}_R}$.
We show that $t_{(j_1,\dots,j_{K-1})}=0$ for all  $(j_1,\dots,j_{K-1})\in S^{K-1}_R$. We fix $(j_1,\dots,j_{K-1})\in S^{K-1}_R$ and set
$\{i_1,\dots, i_{R-K+1}\}=\{1,\dots,R\}\setminus\{j_1,\dots,j_{K-1}\}$. Since $k_{\mathbf C}=K$,  \eqref{eq:2.9} holds. In particular,
\begin{equation}
\alpha_{(j_1,\dots,j_{K-1})}:=\prod\limits_{l=1}^{R-K+1}\det \left[\begin{matrix}\mathbf c_{j_1}&\dots& \mathbf c_{j_{K-1}}&\mathbf c_{i_l}\end{matrix}\right]\ne 0.
\label{eq:2.10}
\end{equation}
Let $\mathbf f:=\mathbf c_{i_1}\otimes\dots\otimes \mathbf c_{i_{R-K+1}}$. Then by \eqref{eq:2.9}, we have
\begin{gather*}
0={\mathbf f}^T\vzero={\mathbf f}^T\left(\mathcal B(\mathbf C)^{(R-K+1)}\widehat{\mathbf t}\right)
=\left(\mathbf c_{i_1}\otimes\dots\otimes \mathbf c_{i_{R-K+1}}\right)^T\Big(\underbrace{\mathcal B(\mathbf C)\odot\dots\odot \mathcal B(\mathbf C)}_{R-K+1}\Big)\widehat{\mathbf t}=\\
\left[\begin{matrix}
\prod\limits_{l=1}^{R-K+1}\det \left[\begin{matrix}\mathbf c_1&\dots& \mathbf c_{K-1}&\mathbf c_{i_l}\end{matrix}\right]
&\dots&
\prod\limits_{l=1}^{R-K+1}\det \left[\begin{matrix}\mathbf c_{R-K+2}&\dots& \mathbf c_{R}&\mathbf c_{i_l}\end{matrix}\right]\end{matrix}\right]{\widehat{\mathbf t}}
=\\
\left[\begin{matrix}
0&\dots&0&
\prod\limits_{l=1}^{R-K+1}\det \left[\begin{matrix}\mathbf c_{j_1}&\dots& \mathbf c_{j_{K-1}}&\mathbf c_{i_l}\end{matrix}\right]&0&\dots&0\end{matrix}\right]{\widehat{\mathbf t}}
=\\
\alpha_{(j_1,\dots,j_{K-1})}\cdot t_{(j_1,\dots,j_{K-1})}.
\end{gather*}
Hence, by \eqref{eq:2.10}, $t_{(j_1,\dots,j_{K-1})}=0$. Thus, $\widehat{\mathbf t}=\vzero$. Therefore, the matrix $\mathcal B(\mathbf C)^{(R-K+1)}$ has full column rank.
\end{proof}

{\em Proof of Proposition 1.10.}
\begin{itemize}
\item[\textup{(i)}] follows from Lemma \ref{lemma:Laplace} \textup{(iv)}.
\item[\textup{(ii)}] follows from Lemma \ref{lemma:2.4} \textup{(iii)}.\qquad\endproof
\end{itemize}

\section{Supplementary material related to properties \textup{(P1)}--\textup{(P4)}}\qquad\\
We will say that condition $\textup{(K,R)}$ holds if for any $K\times R$ matrix $\mathbf C$  with $k_{\mathbf C}=K$ and for any nonzero vector $\mathbf x\in\mathbb R^K$, the
implication in the following scheme holds
\begin{equation}\label{eq:2.88}
\begin{split}
\omega(\mathbf x^T\mathcal B(\mathbf C))&\leq C^{K-1}_R-C^{K-2}_{R-1}=C^{K-1}_{R-1} \Leftrightarrow
\\
&
\mathbf x\ \text{ is orthogonal to at least }\ C^{K-2}_{R-1}\ \text{ columns of }\  B(\mathbf C)\Rightarrow\\
&\mathbf x\ \text{ is proportional to a column of }\ \mathbf C.
\end{split}
\end{equation}
Note that the equivalence ``$\Leftrightarrow$'' and the implication opposite to ``$\Rightarrow$'' in \eqref{eq:2.88} follow from the definition of $\omega(\cdot)$ and Lemma \ref{lemma:2.4} \textup{(ii)}, respectively.
It can be easily checked that $\textup{(2,R)}$ holds for $R\geq 2$ and that
 $\textup{(K,K)}$ holds.
Our goal is to show that $\textup{(K,R)}$ holds for $R\geq K$. We need the following lemma.
\begin{lemma}\label{lemma:2.15}
Suppose that both conditions $\textup{(K-1,R-1)}$ and $\textup{(K,R-1)}$ hold. Then condition $\textup{(K,R)}$ holds.
\end{lemma}
\begin{proof}
Let $\mathbf x\in\mathbb R^{K}$ be a nonzero vector such that
$\omega(\mathbf x^T\mathcal B(\mathbf C))\leq C^{K-1}_{R-1}$ and let
\begin{align*}
\mathcal J&:=\{(j_1,\dots,j_{K-1}):\ (j_1,\dots,j_{K-1})\in S^{K-1}_R,\
\det \left[\begin{matrix}\mathbf c_{j_1}&\dots& \mathbf c_{j_{K-1}}&\mathbf x\end{matrix}\right]=0
\},\\
\mathcal J_1&:=\{(1,j_2,\dots,j_{K-1}):\ (1,j_2,\dots,j_{K-1})\in \mathcal J\},\qquad \overline{\mathcal J}_1:=\mathcal J\setminus\mathcal J_1.
\end{align*}
Then by Lemma \ref{lemma:Laplace} \textup{(i)} and \eqref{eq:2.88}, $\textup{card}\ \mathcal J\geq C^{K-2}_{R-1}$.
We consider two cases:  $\textup{card}\ \overline{\mathcal J}_1\geq C^{K-2}_{R-2}$ and
 $\textup{card}\ \overline{\mathcal J}_1< C^{K-2}_{R-2}$.
If $\textup{card}\ \overline{\mathcal J}_1\geq C^{K-2}_{R-2}$, then
$$
\omega(\mathbf x^T\mathcal B(\left[\begin{matrix}\mathbf c_2&\dots& \mathbf c_R\end{matrix}\right]))\leq C^{K-1}_{R-1}-\textup{card}\ \overline{\mathcal J}_1\leq
C^{K-1}_{R-1}-C^{K-2}_{R-2}=C^{K-1}_{R-2}.
$$
 Since condition $\textup{(K,R-1)}$ holds, it follows that either $\mathbf x$ is zero vector, or
$\mathbf x$ is proportional to one of the vectors $\mathbf c_2,\dots, \mathbf c_R$.
Let us consider the case  $\textup{card}\ \overline{\mathcal J}_1< C^{K-2}_{R-2}$.
We have
\begin{equation}\label{eq:cardJ1}
\textup{card}\ \mathcal J_1=\textup{card}\ \mathcal J-\textup{card}\ \overline{\mathcal J}_1>
 C^{K-2}_{R-1}- C^{K-2}_{R-2}= C^{K-3}_{R-2}.
\end{equation}
By Lemma \ref{lemma:Laplace} \textup{(i)} and \eqref{eq:2.88}, there exist numbers $\alpha_{p_1},\alpha_{p_{j_2}},\dots,\alpha_{p_{j_{K-1}}}$ such that
\begin{equation}\label{eq:T}
\mathbf x=\alpha_{p_1}\mathbf c_1+\sum\limits_{q=2}^{K-1}\alpha_{p_{j_q}}\mathbf c_{j_q},\qquad p\in\{1,\dots,\textup{card}\ \mathcal J_1
\},\qquad (1,j_2,\dots,j_{K-1})\in \mathcal J_1.
\end{equation}
Let $\mathbf T:\ \mathbb R^{K}\rightarrow \mathbb R^{K-1}$ be a linear mapping with $\ker( \mathbf T)=\textup{span}\ \{\mathbf c_1\}$.
We set
$$
\widetilde{\mathbf x}:=\mathbf T\mathbf x\in\mathbb R^{K-1},\qquad
\widetilde{ \mathbf C}:=\mathbf T\left[\begin{matrix}\mathbf c_2&\dots& \mathbf c_R\end{matrix}\right]=
\left[\begin{matrix}\widetilde{\mathbf c}_1&\dots& \widetilde{\mathbf c}_{R-1}\end{matrix}\right]\in\mathbb R^{(K-1)\times (R-1)}.
$$
If  $\mathbf T\mathbf x=\vzero$, then either $\mathbf x=\vzero$, or $\mathbf x$ is proportional to $\mathbf c_1$.
Hence, we can assume that $\widetilde{\mathbf x}$ is a nonzero vector. Since $k_{\mathbf C}=K$ and $\ker(\mathbf T)=\textup{span}\ \{\mathbf c_1\}$, it follows that
$k_{\widetilde{\mathbf C}}=K-1$. Let us apply $\mathbf T$ to \eqref{eq:T}
$$
\widetilde{\mathbf x}=\sum\limits_{q=2}^{K-1}\alpha_{p_{j_q}}\mathbf T\mathbf c_{j_q}
=\sum\limits_{q=2}^{K-1}\alpha_{p_{j_q}}\widetilde{\mathbf c}_{j_q-1}
,\quad p\in\{1,\dots,\textup{card}\ \mathcal J_1
\},\quad (1,j_2,\dots,j_{K-1})\in \mathcal J_1.
$$
Hence,
$$
\det \left[\begin{matrix}\widetilde{\mathbf c}_{j_2-1}&\dots& \widetilde{\mathbf c}_{j_{K-1}-1}&\mathbf x\end{matrix}\right]=0,
\quad (1,j_2,\dots,j_{K-1})\in \mathcal J_1.
$$
By Lemma \ref{lemma:Laplace} \textup{(i)} and \eqref{eq:cardJ1},
$$
\omega(\widetilde{\mathbf x}^T\mathcal B(\widetilde{\mathbf C}))\leq C^{K-2}_{R-1}-\textup{card}\ \mathcal J_1\leq C^{K-2}_{R-1}-C^{K-3}_{R-2}=C^{K-2}_{R-2}.
$$
Since condition $\textup{(K-1,R-1)}$ holds and $\widetilde{\mathbf x}\not=\vzero$, it follows that  $\widetilde{\mathbf x}$ is  proportional to a column of $\widetilde{\mathbf C}$.
Hence, $\mathbf x$ is proportional to one of the vectors $\mathbf c_2,\dots,\mathbf c_R$.
\end{proof}
\begin{lemma}\label{lemma:2.6}
Let $R\geq K$. Then condition $\textup{(K,R)}$ holds.
\end{lemma}

{\em Proof.}
The proof is by induction on $k=2,\dots,K$. For $k=2$ the result is trivial. Suppose that condition $\textup{(k,R)}$ holds for $R\geq k$. We prove that condition
$\textup{(k+1,R)}$ holds for $R\geq k+1$. Since $\textup{(k,k)}$ holds for $k\geq 2$, by Lemma \ref{lemma:2.15}, we have
\begin{align*}
&\text{ conditions }\ \textup{(k,k+1)}\ \text{ and }\ \textup{(k+1,k+1)}\ \text{ imply condition }\      \textup{(k+1,k+2)},\\
&\text{ conditions }\ \textup{(k,k+2)}\ \text{ and }\ \textup{(k+1,k+2)}\ \text{ imply condition }\      \textup{(k+1,k+3)},\\
&\qquad\qquad\qquad\qquad\qquad\qquad\qquad\vdots\\
&\text{ conditions }\ \textup{(k,R-1)}\ \text{ and }\ \textup{(k+1,R-1)}\ \text{ imply condition }\      \textup{(k+1,R)}.\qquad\endproof
\end{align*}
{\em Proof of properties \textup{(P1)}--\textup{(P4)}.}
\begin{itemize}
\item[\textup{(P1)}] follows from Lemma \ref{lemma:2.4} \textup{(i)}.
\item[\textup{(P2)}] follows from Lemma \ref{lemma:Laplace} \textup{(iii)}.
\item[\textup{(P3)}] follows from Lemma \ref{lemma:2.4} \textup{(ii)}.
\item[\textup{(P4)}]
We prove that
\begin{equation*}
\begin{split}
&\mathbf x\ \text{ is orthogonal to }\ C^{K-2}_{R-1}\ \text{ columns of }\  B(\mathbf C)\Leftrightarrow\\
&\mathbf x\ \text{ is proportional to a column of }\ \mathbf C.
\end{split}
\end{equation*}
The result follows from
Lemma \ref{lemma:2.4} \textup{(ii)}  and Lemma  \ref{lemma:2.6}, respectively.
\end{itemize}
\qquad\endproof
\section{Supplementary material related to Lemma 2.17}\qquad\\
{\em Proof of Lemma 2.17.}\\
\textup{(i)}\
From the definitions of the matrices
 $\mathcal Q_m(\mathbf C)^T$ and $\mathcal R_m(\mathbf C)^T$
 (see Definitions 2.9 and 2.10)
 it follows that the matrix $\mathcal Q_m(\mathbf C)^T$ is obtained from $\mathcal R_m(\mathbf C)^T$
be removing columns that are repeated. Thus, \textup{(i)} just expresses the following fact:
$(j_1,\dots,j_m)$-th column of the matrix $\mathcal Q_m(\mathbf C)^T$ coincides with the  $(l_1,\dots,l_m)$-th column of the matrix $\mathcal R_m(\mathbf C)^T$ whenever
$(l_1,\dots,l_m)\in P_{\{j_1,\dots,j_m\}}$.

\textup{(ii)}\ From Lemma 2.15 and definition of the matrix $\mathbf G$ (2.5) it follows that
the $\left((i_1,\dots,i_m),(j_1,\dots,j_m)\right)$-th entry of the matrix $\mathcal R_m(\mathbf C)^T\mathbf G$ is equal to
\begin{align*}
&m!\left(\pi_S(\mathbf c_{i_1}\otimes\dots\otimes\mathbf c_{i_m})\right)^T\left(\pi_S(\mathbf e^K_{j_1}\otimes\dots\otimes\mathbf e^K_{j_m})\right)=\\
&\frac{1}{m!}
\sum\limits_{(q_1,\dots,q_m)\in P_{\{i_1,\dots,i_m\}}}
\sum\limits_{(l_1,\dots,l_m)\in P_{\{j_1,\dots,j_m\}}}
(\mathbf c_{q_1}^T\mathbf e_{l_1}^K)\cdot\dots\cdot (\mathbf c_{q_m}^T\mathbf e_{l_m}^K)=\\
&\frac{1}{m!}
\sum\limits_{(q_1,\dots,q_m)\in P_{\{i_1,\dots,i_m\}}}
\textup{perm\ }\mathbf C((j_1,\dots,j_m),(q_1,\dots,q_m)).
\end{align*}
Since $\textup{perm\ }\mathbf C((j_1,\dots,j_m),(q_1,\dots,q_m))=\textup{perm\ }\mathbf C((j_1,\dots,j_m),(i_1,\dots,i_m))$ for all
$(q_1,\dots,q_m)\in P_{\{i_1,\dots,i_m\}}$, it follows that
\begin{equation*}
\left(\pi_S(\mathbf c_{i_1}\otimes\dots\otimes\mathbf c_{i_m})\right)^T\left(\pi_S(\mathbf e^K_{j_1}\otimes\dots\otimes\mathbf e^K_{j_m})\right)=
\textup{perm\ }\mathbf C((j_1,\dots,j_m),(i_1,\dots,i_m)).
\end{equation*}
The equality in \textup{(ii)} follows now from Definition 2.9.
\qquad\endproof
\section{Supplementary material related to Corollary 2.18}\qquad\\
{\em Proof of Corollary 2.18.}\\
By Proposition 2.13 \textup{(ii)}, $\dim \left( \ker \left(\mathcal{R}_m(\mathbf C)^T\upharpoonright_{\textup{range}(\pi_S)}\right)\right)=C^{K-1}_R$. We prove that 
the equation $\mathcal{R}_m(\mathbf C)^T\mathbf x=\vzero$ has
at most $C^{K-1}_R-(K-1)$ solutions of the form
$\mathbf x:=\mathbf y\otimes\mathbf z\in\textup{range}(\pi_S)$.
 Denote by $\mathcal X$  the  $m$-th order symmetric  
$K\times\dots\times K$ tensor whose vectorized version coincides with $\mathbf x$.
Since $\mathcal X$ is symmetric, all its $K\times K^{m-1}$ matrix unfoldings
coincide with the rank-1 matrix $\mathbf y\mathbf z^T$ and hence are rank-1.
It is well known that this is possible if and only if  $\mathcal X$  is itself a rank-1 tensor.
Hence, the vector $\mathbf x$ is proportional to the vector $\mathbf y\otimes\dots\otimes\mathbf y$ and
$\mathcal{R}_m(\mathbf C)^T\mathbf x=\mathcal R_m(\mathbf C)^T(\mathbf y\otimes\dots\otimes \mathbf y)=\vzero$.
Hence, by (2.9), $(\mathbf c_{p_1}^T\mathbf y)\cdots (\mathbf c_{p_m}^T\mathbf y)=0$ for
$1\leq p_1<\dots<p_m\leq R$.
Consequently, $\mathbf y$ is orthogonal to at least $R-(m-1)=K-1$ columns of $\mathbf C$.
Since $k_{\mathbf C}=K-1$, the orthogonal complement of these $K-1$ columns is one-dimensional and, hence, 
by Lemma \ref{lemma:Laplace} \textup{(iii)}, $\mathbf y$ is proportional to a column of $\mathcal B(\mathbf C)$.
On the other hand, since $k_{\mathbf C}= K-1$, $\mathbf C$ has $K$ columns that are linearly dependent, and, hence,
 it follows from Lemma \ref{lemma:Laplace} \textup{(iii)} that at least $K$ columns of $\mathcal B(\mathbf C)$ are pairwise proportional.
Hence, there exist at most $C^{K-1}_R-(K-1)$ nonzero solutions of the equation $\mathcal R_m(\mathbf C)^T(\mathbf y\otimes\dots\otimes \mathbf y)=\vzero$, which completes the proof.
\qquad\endproof
\section{Supplementary material related to Lemma 4.4}\qquad\\
We need the following Lemma.
\begin{lemma}\label{lemma:before:4.3}
Let $\widehat{\mathbf A}$ be any set of $\widehat{m}$ columns of $\mathbf A$, let
$\widehat{\mathbf B}$ be the corresponding set of columns of $\mathbf B$.
Assume that the matrix $\mathcal C_m(\mathbf A)\odot \mathcal C_m(\mathbf B)$  has full column rank.
Then
\begin{itemize}
\item[\textup{(i)}] $\min(k_{\mathbf A},k_{\mathbf B})\geq m$.
\item[\textup{(ii)}] $\max(r_{\widehat{\mathbf A}},r_{\widehat{\mathbf B}})\geq\min(\widehat{m},m+1)$.
\end{itemize}
\end{lemma}
\begin{proof}
Since $\mathcal C_m(\mathbf A)\odot \mathcal C_m(\mathbf B)$  has full column rank it follows that
all columns of $\mathcal C_m(\mathbf A)$ and $\mathcal C_m(\mathbf B)$ are nonzero. Hence,
\textup{(i)} follows from Lemma 2.3 (1). If $\widehat{m}\leq m$, then
\textup{(ii)} follows from \textup{(i)}. If $\widehat{m}> m$, then, by
Lemma 2.5, $r_{\mathcal C_m(\widehat{\mathbf A})}=C^m_{r_{\widehat{\mathbf A}}}$
and $r_{\mathcal C_m(\widehat{\mathbf B})}=C^m_{r_{\widehat{\mathbf B}}}$. Hence,
\begin{equation}\label{eq:4.5new}
\begin{split}
C^m_{r_{\widehat{\mathbf A}}}C^m_{r_{\widehat{\mathbf B}}}=
r_{\mathcal C_m(\widehat{\mathbf A})}r_{\mathcal C_m(\widehat{\mathbf B})}=
r_{\mathcal C_m(\widehat{\mathbf A})\otimes \mathcal C_m(\widehat{\mathbf B})}\geq
r_{\mathcal C_m(\widehat{\mathbf A})\odot \mathcal C_m(\widehat{\mathbf B})}=
C^m_{\widehat{m}}>1,
\end{split}
\end{equation}
where the last equality holds since the matrix $\mathcal C_m(\widehat{\mathbf A})\odot \mathcal C_m(\widehat{\mathbf B})$
has full column rank. The statement \textup{(ii)} for $\widehat{m}> m$
now follows from \eqref{eq:4.5new}.
\end{proof}

{\em Proof of Lemma 4.4.}\\
 Without loss of generality we may assume that $\mathbf F$ coincides with $\mathcal B(\mathbf C)$. In the proof we will associate
 indices $i,j\in\{1,\dots,C^{K-1}_R\}$ with their multi-index analogues
$( i_1,\dots,i_{K-1}), ( j_1,\dots,j_{K-1})\in S^{K-1}_R$.
 By definition set
 $$
 \{p_1,\dots,p_{\widehat{m}}\}=\{1,\dots,R\}
 \setminus\left\{\{i_1,\dots,i_{K-1}\}\cap\{j_1,\dots,j_{K-1}\}\right\}.
 $$
  Since
 \begin{gather*}  \textup{card}\left\{\{i_1,\dots,i_{K-1}\}\cap\{j_1,\dots,j_{K-1}\}\right\}\leq K-2,\\
 \{p_1,\dots,p_{\widehat{m}}\}
 \subset\left\{\{1,\dots,R\}\setminus\{i_1,\dots,i_{K-1}\}\right\}
 \cup
 \left\{\{1,\dots,R\}\setminus\{j_1,\dots,j_{K-1}\}\right\}
\end{gather*}
 it follows that
 \begin{gather}
 m=R-(K-2)\leq \widehat{m},\label{eqmhat1}\\
\widehat{m}\leq \min(R, 2(R-(K-1)))=\min(R,2m-2).
 \label{eqmhat2}
 \end{gather}
  We will show that the statements \textup{(i)}--\textup{(iii)} are all equivalent to the condition $\widehat{m}=m$.

 \textup{(i)} $\Leftrightarrow \widehat{m}=m$:\
 follows from Lemma \ref{lemma:2.4} \text{(i)}.

  \textup{(ii)} $\Leftrightarrow \widehat{m}=m$:\
let the vectors $\widehat{\mathbf y}_i$ and $\widehat{\mathbf y}_j$ (resp. the matrices
 $\widehat{\mathbf A}$ and $\widehat{\mathbf B}$) be formed by
the entries of the vectors  ${\mathbf y}_i$ and ${\mathbf y}_j$ (resp. by the   columns of the matrices  $\mathbf A$ and $\mathbf B$) with indices $p_1,\dots,p_{\widehat{m}}$.
Then
\begin{equation}
\mathbf V_i=\mathbf A\Diag{{\mathbf y}_i}\mathbf B^T=\widehat{\mathbf A}
\Diag{{\widehat{\mathbf y}}_i}
\widehat{\mathbf B}^T,\ \
\mathbf V_j=\mathbf A\Diag{{\mathbf y}_j}\mathbf B^T=
\widehat{\mathbf A}\Diag{{\widehat{\mathbf y}}_j}
\widehat{\mathbf B}^T.\label{eq:4.5}
\end{equation}
Hence,
 \begin{gather}
[\mathbf V_i\ \mathbf V_j]=
\widehat{\mathbf A}[
\Diag{{\widehat{\mathbf y}}_i
}\widehat{\mathbf B}^T\
\Diag{{\widehat{\mathbf y}}_j
}\widehat{\mathbf B}^T]=
\widehat{\mathbf A}(\widehat{\mathbf C}\odot \widehat{\mathbf B})^T,
\label{eq:mr1}\\
[\mathbf V_i^T\ \mathbf V_j^T]=
\widehat{\mathbf B}[
\Diag{{\widehat{\mathbf y}}_i
}\widehat{\mathbf A}^T\
\Diag{{\widehat{\mathbf y}}_j
}\widehat{\mathbf A}^T]=
\widehat{\mathbf B}(\widehat{\mathbf C}\odot \widehat{\mathbf A})^T,
\label{eq:mr2}
\end{gather}
where $\widehat{\mathbf C}:=[\widehat{\mathbf y}_i\ \widehat{\mathbf y}_j]^T$.
We claim  that
\begin{equation}
 \text{the matrices }\
\widehat{\mathbf C}\odot \widehat{\mathbf B}\ \text{ and }\
\widehat{\mathbf C}\odot \widehat{\mathbf A}\ \text{ have full column rank.}\label{eq:CkrBandCkrAfcr}
\end{equation}
From the construction of the matrix $\widehat{\mathbf C}$ it follows that there exists an $\widehat{m}\times\widehat{m}$ permutation matrix $\mathbf P$ such that
$$
(\widehat{\mathbf C}\odot \widehat{\mathbf B}){\mathbf P}=
(\widehat{\mathbf C}\mathbf P)\odot (\widehat{\mathbf B}{\mathbf P})=
\left[
\begin{matrix}
*&\dots&*&*&\dots&*&0&\dots&0\\
\coolunder{_{\widehat{m}-m+1}}{0&\dots&0}&\coolunder{_{2m-\widehat{m}-2}}{*&\dots&*}&\coolunder{_{\widehat{m}-m+1}}{*&\dots&*}
\end{matrix}
\right]\odot (\widehat{\mathbf B}{\mathbf P})=
$$
$$
\left[
\begin{matrix}
*\cdot\widetilde{\mathbf b}_1&\dots&*\cdot\widetilde{\mathbf b}_{\widehat{m}-m+1}&*\cdot\widetilde{\mathbf b}_{\widehat{m}-m+2}&\dots&*\cdot\widetilde{\mathbf b}_{m-1}&0&\dots&0\\
0        &\dots&0                &*\cdot\widetilde{\mathbf b}_{\widehat{m}-m+2}&\dots&*\cdot\widetilde{\mathbf b}_{m-1}&*\cdot\widetilde{\mathbf b}_{m}&\dots&*\cdot\widetilde{\mathbf b}_{\widehat{m}}
\end{matrix}
\right],
$$
where $*$ denotes a nonzero value, $\widetilde{\mathbf b}_1,\dots,\widetilde{\mathbf b}_{\widetilde{m}}$ denote the  columns of the
matrix $\widehat{\mathbf B}{\mathbf P}$, and where we use the dimensionality constraints  in \eqref{eqmhat1}--\eqref{eqmhat2}.
By Lemma \ref{lemma:before:4.3} \textup{(i)}, $k_{\widehat{\mathbf B}}\geq k_{\mathbf B}\geq m$. Hence, the matrix $(\widehat{\mathbf C}\odot \widehat{\mathbf B})\mathbf P$ has full column rank. Since the matrix
$\mathbf P$ is nonsingular, it follows that the matrix $\widehat{\mathbf C}\odot \widehat{\mathbf B}$ also has full column rank.
In a similar fashion one can prove that
the matrix $\widehat{\mathbf C}\odot \widehat{\mathbf A}$  has full column rank.
From \eqref{eq:mr1}--\eqref{eq:CkrBandCkrAfcr} it follows that
$
r_{[\mathbf V_i\ \mathbf V_j]}=r_{\widehat{\mathbf A}}$ and
$r_{[\mathbf V_i^T\ \mathbf V_j^T]}=r_{\widehat{\mathbf B}}$.
Using  Lemma \ref{lemma:before:4.3} we have
\begin{gather}
\widehat{m}\geq\max(r_{[\mathbf V_i\ \mathbf V_j]},r_{[\mathbf V_i^T\ \mathbf V_j^T]})=
\max(r_{\widehat{\mathbf A}},r_{\widehat{\mathbf B}})\geq\min(\widehat{m},m+1),\label{eq:munfmax}\\
\min(r_{[\mathbf V_i\ \mathbf V_j]},r_{[\mathbf V_i^T\ \mathbf V_j^T]})=
\min(r_{\widehat{\mathbf A}},r_{\widehat{\mathbf B}})
\geq\min(k_{\widehat{\mathbf A}},k_{\widehat{\mathbf B}})
\geq\min(k_{\mathbf A},k_{\mathbf B})
\geq m.\label{eq:munfmin}
\end{gather}
The equivalence \textup{(ii)} $\Leftrightarrow \widehat{m}=m$ now easily follows from \eqref{eq:munfmax}--\eqref{eq:munfmin}.

\textup{(iii)} $\Leftrightarrow \widehat{m}=m$:\ by \eqref{eq:4.5},
$\mathcal V_{ij}=
[\mathbf A,\mathbf B, [{\mathbf y}_i\ {\mathbf y}_j]^T]_R=
[\widehat{\mathbf A},\widehat{\mathbf B},\widehat{\mathbf C}]_{\widehat{m}}$. Hence, $\widehat{m}\geq r_{\mathcal V_{ij}}$. On the other hand, from \eqref{eq:mr1}-\eqref{eq:mr2}
it follows that  $r_{\mathcal V_{ij}}\geq \
\max(r_{[\mathbf V_i\ \mathbf V_j]},$ $r_{[\mathbf V_i^T\ \mathbf V_j^T]})$.
Hence, by \eqref{eq:munfmax}, $r_{\mathcal V_{ij}}\geq \min(\widehat{m},m+1)$. The equivalence \textup{(iii)} $\Leftrightarrow \widehat{m}=m$ now follows from the inequalities
$\widehat{m}\geq r_{\mathcal V_{ij}}\geq \min(\widehat{m},m+1)$.
\qquad\endproof
\section{Supplementary material related to Example 4.6}\qquad\\
The full  matrix $\TTm{3}$ is given by
$$
\TTm{3}=-
\left[
\begin{array}{r@{\hspace{2pt}}r@{\hspace{2ex}}r@{\hspace{2ex}}r@{\hspace{2pt}}r@{\hspace{2pt}}r@{\hspace{2pt}}r@{\hspace{2ex}}r@{\hspace{2ex}}r@{\hspace{2ex}}r@{\hspace{2ex}}r@{\hspace{2pt}}r@{\hspace{2pt}}r@{\hspace{2pt}}r@{\hspace{2pt}}r@{\hspace{2pt}}r@{\hspace{2ex}}r@{\hspace{2pt}}r@{\hspace{2pt}}r@{\hspace{2ex}}r}
0&  0&	 0&	0&	0&  {\bf 0}&    6&	0&	0&	0& 0&   0&  12&	0&  6& 12& 0&  0&  0& 0\\
0&  0&	 0&	0&	0&  {\bf 0}&	0&	0&  2&	0& 0&	0&   0&	0& -2&	0& 0& -4& -4&0\\
0&  0&	 0&	0&	0&  {\bf 3}&	0&	0&	0&	0& 0&  -6&   0&-6&  3&	0& 0&  0&  0&0\\
0&  0&	 0&	0&	0&  {\bf 3}&	0&  0&	2&	0& 0&	6&	 0&	6& 11&	0& 0&  4&  4&0\\
0& -6&	 0&	0& -6&  {\bf-3}&   -3&  0&	0&	0& 0&	0&   0&	0&  0&	0& 0&  0&	 0&0\\
0&  0&	 2&	0&	0&  {\bf 1}&	0&	2&	1&	0& 0&	0&	 0&	0&  0&	0& 0&  0&	 0&0\\
0&  0&	 2&	0&	0&  {\bf 4}&	0&	2&	1&	0& 0&	0&	 0&	0&  0&  0& 0&  0&	 0&0\\
0&  0&   0& 0&  0&  {\bf 0}&    0&  0&  0&  0& 0&   0&   0& 0&  0&  0& 0&  0&     0&0\\
0&  0&	 0&	4&	0&  {\bf 0}&	2&	0&	2&	4& 0&	0&   0&	0&  0&	0& 0&  0&	 0&0\\
0&  0&	 2&	0&	0&  {\bf 1}&	0&	2&	3&	0& 0&	0&	 0&	0&  0&	0& 0&  0&	 0&0\\
0&  0&	 2& 0&	0&  {\bf 1}&	0&	2&	1&	0& 0&	0&	 0&	0&  0&	0& 0&  0&	 0&0\\
0&  0&   0& 0&  0&  {\bf 0}&    0&  0&  0&  0& 0&   0&   0& 0&  0&  0& 0&  0&     0& 0\\
0&  6&	 0&	4&	6&  {\bf 3}&   11&  0&  2&  4& 0&	0&	 0&	0&  0&	0& 0&  0&	0& 0\\
0&  0&   0& 0&  0&  {\bf 0}&    0&  0&  0&  0& 0&   0&   0& 0&  0&  0& 0&  0&     0&0\\
0&  0&   0& 0&  0&  {\bf 0}&    0&  0&  0&  0& 0&   0&   0& 0&  0&  0& 0&  0&     0& 0\\
0&  0&   0& 0&  0&  {\bf 0}&    0&  0&  0&  0& 0&   0&   0& 0&  0&  0& 0&  0&     0& 0
\end{array}
\right].
$$